\documentclass[12pt,a4paper]{amsart}
\usepackage{geometry}
\usepackage{tabls}

\usepackage{url}

\usepackage{amssymb}
\usepackage{mathrsfs}
\usepackage{amsmath}
\usepackage{amsthm}
\usepackage{amscd}
\usepackage{fancyhdr}
\usepackage{enumerate}
\usepackage{paralist}
\usepackage{graphicx}
\usepackage{cite}

\usepackage{footmisc}
\numberwithin{equation}{section}

\theoremstyle{plain}
\newtheorem{Cor}{Corollary}[section]

\newtheorem{Thm}[equation]{Theorem}
\newtheorem{lem}[equation]{Lemma}
\newtheorem{prop}[equation]{Proposition}
\newtheorem{rem}[equation]{Remark}

\begin{document}

\title{The topology of the set of multiple zeta-star values}

\author{Jiangtao Li}

\email{lijiangtao@csu.edu.cn}
\address{Jiangtao Li \\ School of Mathematics and Statistics, HNP-LAMA, Central South University, Hunan Province, China}

\begin{abstract}
  We provide a multiple integral representation for each multiple zeta-star value, and utilize these representations to establish a natural order structure on the set of such values. This order structure allows for a one-to-one correspondence between a subset of the infinite sequences of natural numbers and the half line $(1,+\infty)$.  Some basic properties of this correspondence are discussed. We also calculate the Hausdorff dimensions for the images of some subsets of the  infinite sequences under this correspondence. As a result of this correspondence, we are able to determine the limits for a number of natural multiple integrals.  Our analysis also reveals that the set of multiple zeta-star values is dense within the $(1,+\infty)$ domain, and that each value is non-integer in nature. 
 \end{abstract}

\let\thefootnote\relax\footnotetext{
2020 $\mathnormal{Mathematics} \;\mathnormal{Subject}\;\mathnormal{Classification}$. 11M32, 11K55. \\
$\mathnormal{Keywords:}$  Multiple zeta-star values,  Hausdorff dimension.\\
Project funded by the National Natural Science Foundation of China (Grant No. 12571009). }

\maketitle

\section{Introduction}\label{int}
      Multiple zeta values are defined by
      \[
      \zeta(k_1,\cdots,k_r)=\sum_{n_1>\cdots>n_r>0}\frac{1}{n_1^{k_1}\cdots n_r^{k_r}}, k_1\geq 2,k_2,\cdots, k_{r}\geq 1.
      \]
      For multiple zeta value $\zeta(k_1,\cdots,k_r)$, denote by $N=k_1+\cdots+k_r$ and $r$ its weight  and depth respectively.

  For  $k_1\geq 2,k_2,\cdots, k_{r}\geq 1$, the multiple zeta-star value $\zeta^{\star}(k_1,\cdots,k_r)$ is defined by
        \[
        \zeta^{\star}(k_1,\cdots,k_r)=\sum_{n_1\geq\cdots\geq n_r\geq 1}\frac{1}{n_1^{k_1}\cdots n_r^{k_r}}.
        \]
        It is easy to check that multiple zeta-star values are $\mathbb{Q}$-linear combinations of multiple zeta values. Every multiple zeta value is also a $\mathbb{Q}$-linear combination of multiple zeta-star values. 
        
        Every multiple zeta value has an iterated path integral representation. In this paper, we show that every multiple zeta-star value also has  a multiple integral representation.

 \begin{Thm}\label{last}
(i) For $r=2k+1$, $1\leq i_1<\cdots<i_r$,
\[
\begin{split}
&\;\;\mathop{\int}_{[0,1]^{i_r}}\frac{dx_1\cdots dx_{i_r}}{1-x_1\cdots x_{i_1}+x_1\cdots x_{i_2}+\cdots+(-1)^rx_1\cdots x_{i_r}}\\
&= \zeta^{\star}( i_1+1   , \underbrace{1,\cdots,1}_{i_2-i_1-1},    i_3-i_2+1, \underbrace{1,\cdots, 1}_{i_4-i_3-1},  \cdots,   i_{2k-1}-i_{2k-2}+1,  \underbrace{1,\cdots,1}_{i_{2k}-i_{2k-1}-1},  i_{2k+1}-i_{2k}      ),
\end{split}
\]
(ii) For $r=2k$, $1\leq i_1<\cdots<i_r$,
\[
\begin{split}
&\;\;\mathop{\int}_{[0,1]^{i_r}}\frac{dx_1\cdots dx_{i_r}}{1-x_1\cdots x_{i_1}+x_1\cdots x_{i_2}+\cdots+(-1)^rx_1\cdots x_{i_r}}\\
&=\zeta^{\star}( i_1+1   , \underbrace{1,\cdots,1}_{i_2-i_1-1},    i_3-i_2+1, \underbrace{1,\cdots, 1}_{i_4-i_3-1},  \cdots,   i_{2k-1}-i_{2k-2}+1,  \underbrace{1,\cdots,1}_{i_{2k}-i_{2k-1}-1}).
\end{split}
\]
\end{Thm}

In fact, the above formula was already known to Zlobin\cite{zlo}. 

For a subset $\mathcal{C}\subseteq \mathbb{R}$, denote by $\mathcal{C}^\prime$ the set of accumulation points of $\mathcal{C}$. Inductively, define 
\[
\mathcal{C}^{(1)}=  \mathcal{C}^\prime     , \mathcal{C}^{(n+1)}=\left( \mathcal{C}^{(n)}  \right)^\prime, \;r\geq 1.
\]
The subset $C^{(r)}$ is called the $n$-th derived set of $\mathcal{C}$.
 
Kumar \cite{kum} studied the order structure and the $n$-th derived sets for the set of multiple zeta values.
Based on Theorem \ref{last}, for the order structure of the set of multiple zeta-star values, we have:
\begin{Thm}\label{order}
Denote by 
\[
\mathcal{S}=\Big{\{}(k_1,\cdots,k_r)\,\Big{|}\,k_1\geq 2, k_2,\cdots,k_r\geq 1, r\geq1\Big{\}}.
\]
Define an order $\succ$ on $\mathcal{S}$ by
\[
(k_1,\cdots, k_{r},k_{r+1})\succ(k_1,\cdots, k_{r}), \forall\, (k_1,\cdots,k_r, k_{r+1})\in \mathcal{S}\]
and 
\[
(k_1,\cdots, k_{r})\succ(m_1,\cdots, m_{s})
\]
if $  k_i=m_i, 1\leq i\leq j, k_{j+1}<m_{j+1}$, for some $j\geq 0$.
 Then for any $$(k_1, \cdots, k_r), (m_1, \cdots, m_s)\in \mathcal{S},$$
 \[
 (k_1,\cdots, k_{r})\succ(m_1,\cdots, m_{s}) \]
 if and only if 
 \[
  \zeta^{\star}(k_1,\cdots, k_{r})>\zeta^{\star}(m_1,\cdots, m_{s}). 
  \]
\end{Thm}
 
 Theorem \ref{order} shows that for any $x\in \mathbb{R}$, there is at most ${\bf one}$ multiple zeta-star value $ \zeta^{\star}(k_1,\cdots, k_{r})$ which satisfies $ \zeta^{\star}(k_1,\cdots, k_{r})=x$. Besides, if  $$\zeta^{\star}(k_1,\cdots, k_{r})=\zeta^{\star}(m_1,\cdots, m_{s}),$$ then $r=s$ and $(k_1,\cdots, k_{r})=(m_1,\cdots, m_{s})$. So the formulas of the form 
 \[
 \zeta(2,1)=\zeta(3)
 \]
 don't exist in the set of multiple zeta-star values.
 
 Denote by $\mathcal{Z}^{\star}$ the set of multiple zeta-star values. It is clear that $$\mathcal{Z}^{\star}\subseteq (1,+\infty).$$
 Define 
 \[
 \mathcal{T}=\Big{\{} (k_1,k_2, \cdots, k_r,\cdots)\,\Big{|}\,k_1\geq 2,k_i\geq 1, i\geq 2, k_s\geq 2 \;for \; some \; s\geq 2 \;if \;k_1=2  , k_i\in\mathbb{N}    \Big{\}}.
 \]
 \begin{Thm}\label{rep} 
$(i)$ For any ${\bf k}=(k_1,k_2, \cdots, k_r,\cdots)\in  \mathcal{T}$, the limit 
  \[
 I=\mathop{\mathrm{lim}}_{r\rightarrow +\infty}\zeta^{\star}(k_1,k_2,\cdots,k_r).
 \]
 exists and  $I>1$.\\
 $(ii)$ Furthermore, the map 
 \[
 \eta: \mathcal{T}\rightarrow (1, +\infty),
  \]
  \[
  {\bf k}=(k_1,k_2, \cdots, k_r,\cdots) \mapsto x=\mathop{\mathrm{lim}}_{r\rightarrow +\infty}\zeta^{\star}(k_1,k_2,\cdots,k_r)  \]
 is bijective. One can extend the order on $\mathcal{S}$ naturally to an unique order on $\mathcal{T}$. Then 
 \[
{\bf k} \succ {\bf m}
 \]
 if and only if 
 \[
 \eta({ \bf k}) >\eta ( {\bf m}). 
 \]
 $(iii)$ As a result, $\left(  \mathcal{Z}^{\star}   \right)^{(n)}= [1, +\infty], \forall\, n\geq 1$.\\
  \end{Thm}
 Theorem \ref{rep} shows that the set $\mathcal{Z}^{\star}$ is dense in $(1, +\infty)$. One can compare Theorem \ref{rep} with the theory of continued fractions for any $u\in (0,1)$.   In \cite{hmo}, M. Hirose, H. Murahara and T. Onozuka also found the bijective map $\eta$ independently. They also calculate the images of many special sequences of non-negative integers.
 
 For $x\in (1,+\infty)$, we are curious about its inverse image $\eta^{-1}(x)$. For a sequence $(k_1,\cdots,k_r,\cdots)\in \mathcal{T}$, if 
 \[
 \mathop{sup}_{r\geq 1}\; k_r=+\infty, 
 \]
 then the sequence $(k_1,\cdots,k_r,\cdots)\in \mathcal{T}$ is called unbounded.
  \begin{Thm}\label{zero}  Denote by $m$ the Lebesgue measure on $\mathbb{R}$.\\
  $(i)$ For $p\geq 2$, define 
 \[
 \mathcal{T}_p= \{(k_1,\cdots,k_r,\cdots)\in\mathcal{T}\,|\, k_i\geq p, \forall \, i\geq 1\}.
 \]
 then
  \[
 m(\eta( \mathcal{T}_p))=0, \forall\, p\geq 2. \]
 $(ii)$ For 
 \[
 \mathcal{E}=\{(k_1,\cdots,k_r,\cdots)\in\mathcal{T}\,|\,k_i\geq 2, \forall \,i\geq j, for\;some\;j\geq 1\},
 \]
 then $\left(\eta(\mathcal{E})\right)^\prime=[1,+\infty]$ and  $m(\eta(\mathcal{E}))=0$. In conclusion, for almost all $x\in (1,+\infty)$, its corresponding sequence $\eta^{-1}(x)$ has infinite  1.\\
 $(iii)$  For $p\geq 2$, define 
 \[
 \mathcal{D}_p= \{(k_1,\cdots,k_r,\cdots)\in\mathcal{T}\,|\, k_i\leq  p, \forall \, i\geq 1\}.
 \]
 then  $\eta(\mathcal{D}_p)$ is closed in $ (1,+\infty)$ and 
     \[
 m(\eta(\bigcup_{p\geq 2} \mathcal{D}_p))=0. \] 
 As a result, for almost all $x\in (1,+\infty)$, its corresponding sequence $\eta^{-1}(x)$ is unbouned. \end{Thm}
 
 By Theorem \ref{zero} $(i)$, the set  $\eta(\mathcal{T}_p)$ is very similar to the Cantor set on $[0,1]$.  Theorem \ref{zero} $(ii), (iii)$ shows that for almost all $x\in (1,+\infty)$, its correspondence 
 \[
 \eta^{-1}(x)=(k_1,\cdots,k_r,\cdots)
 \]
 satisfies 
 \[
  \varliminf_{r\rightarrow +\infty}k_r=1, \varlimsup_{r\rightarrow +\infty}k_r=+\infty.
 \]
 
 It is well-known that the Hausdorff dimension of the Cantor set is $\frac{\mathrm{log}\,2}{\mathrm{log}\,3}$. For the set $\eta(\mathcal{T}_p)$ and $\eta(\mathcal{D}_p)$, we have
 \begin{Thm}\label{hdim}
 For $C\subseteq \mathbb{R}^n$, denote by $\mathrm{dim}_H C$ the Housdorff dimension of $C$.\\
 $(i)$ For $p\geq 2$, define $\alpha_p      $  as the unique root of the equation 
 \[
 x^{p-1}(x-1)=1
 \] which satisfies $\alpha_p\in (1,2)$,
  then
 \[
 \mathrm{dim}_H\, \eta\left(\mathcal{T}_p\right)= \frac{\mathrm{log}\;\alpha_p}{\mathrm{log}\;2}, \forall\,p\geq 2.
 \]
 $(ii)$ For $2\leq p< q$, 
 define $\gamma_{p,q}      $  as the unique root of the equation 
 \[
 x^p+x^{p+1}+\cdots+x^q=1
 \] which satisfies $\gamma_{p,q}\in (0,1)$, then 
  \[
 \mathrm{dim}_H\, \eta\left(\mathcal{T}_p \cap \mathcal{D}_q\right)=\frac{\mathrm{log}\;\frac{1}{\gamma_{p,q}}}{\mathrm{log}\;2}.
 \]
 $(iii)$ For $q\geq 3$, one has 
 \[
  \mathrm{dim}_H\, \eta\left( \mathcal{D}_q\right)\geq \frac{\mathrm{log}\;\frac{1}{\gamma_{2,q}}}{\mathrm{log}\;2}.
 \]
  \end{Thm}

 Gorodetsky, Lichtman and Wong \cite{GLW} showed that
 \[
 \mathop{\mathrm{lim}}_{r\rightarrow +\infty} I_r= \mathop{\mathrm{lim}}_{r\rightarrow +\infty} \mathop{\int}_{[0,1]^r}\frac{dx_1dx_2\cdots dx_r}{1+x_1+\cdots+x_1x_2\cdots x_r    }=e^{-\gamma}.
 \]
 Actually, the integral $I_r$ is equal to the continuous version multiple zeta values $\zeta^{\mathcal{C}}(\{1\}^r,2)$ which was defined by the author in \cite{li} (Proposition $3.4$ in \cite{li}). 
 
 By Theorem \ref{last} and Theorem \ref{rep}, one has the following result.
 \begin{Cor}\label{alt}
 Denote by 
 \[
 \mathcal{A}=\Big{\{} (i_1,i_2,\cdots, i_r,\cdots)\,\Big{|}\, 1\leq i_1<i_2<\cdots<i_r<\cdots, \,i_j\in\mathbb{N}, \forall\, j\geq 1     \Big{\}}
 \]
 $(i)$ For any ${\bf i}=(i_1,i_2,\cdots, i_r,\cdots)\in\mathcal{A}$, the limit 
 \[
I_{\bf i}= \mathop{\mathrm{lim}}_{r\rightarrow +\infty} \mathop{\int}_{[0,1]^{i_r}}\frac{dx_1dx_2\cdots dx_{i_r}}{1-x_1x_2\cdots x_{i_1}+ x_1x_2\cdots x_{i_2}  +\cdots+(-1)^rx_1x_2\cdots x_{i_r}    } \]
 exists and $I_i>1$.\\
 $(ii)$ The map \[\xi: \mathcal{A}\rightarrow (1,+\infty) ,\]
 \[
 (i_1,i_2,\cdots, i_r,\cdots)\mapsto \mathop{\mathrm{lim}}_{r\rightarrow +\infty} \mathop{\int}_{[0,1]^{i_r}}\frac{dx_1dx_2\cdots dx_{i_r}}{1-x_1x_2\cdots x_{i_1}+ x_1x_2\cdots x_{i_2}  +\cdots+(-1)^rx_1x_2\cdots x_{i_r}    }  \]
 is injective. What is more, $\xi(\mathcal{A})=(1,+\infty)-\mathcal{Z}^{\star}$.
 \end{Cor}
 
 \begin{Cor}
 By Theorem \ref{alt},  if $\mathbb{Q}\cap (1,+\infty)\subseteq \xi(\mathcal{A})$, then every multiple zeta-star value is irrational.
 \end{Cor}
 
 In Section \ref{sie}, we will give another approach to Theorem \ref{last} by Yamamoto's result \cite{yam}. In the last section, we will calculate the limits
 \[
 \mathop{\mathrm{lim}}_{r\rightarrow +\infty}\zeta^{\star}(k_1,k_2,\cdots,k_r) \]
 and 
 \[
 \mathop{\mathrm{lim}}_{r\rightarrow +\infty} \mathop{\int}_{[0,1]^{i_r}}\frac{dx_1dx_2\cdots dx_{i_r}}{1-x_1x_2\cdots x_{i_1}+ x_1x_2\cdots x_{i_2}  +\cdots+(-1)^rx_1x_2\cdots x_{i_r}    }  \]
 for some special ${\bf k}\in\mathcal{T}$ and ${\bf i}\in \mathcal{A}$.
 As an application, we show 
 that
 \begin{Thm}\label{itg}
 Every multiple zeta-star value is not an integer.
 \end{Thm}
 Although the set of multiple zeta-star values is dense in $(1,+\infty)$, every multiple zeta-star value is still  not an integer. Theorem \ref{itg} is compatible with the conjecture of transcendence of multiple zeta values.
 
 \begin{rem}
 After this paper was announced at arXiv, there are many further developments on the theory of multiple zeta-star values. See \cite{lir} for the rational deformations of multiple zeta-star values and $n$-th Cantor set. See \cite{lid} for the approximation theory of multiple zeta-star values. See Kamano \cite{kam} for the order structures of multi-polylogarithms.
 \end{rem}

        \section{The order structure of the set of the multiple zeta-star values}\label{gea}
In this section, we will give the multiple integral representation for every multiple zeta-star value. Based on these multiple integral representations, we will characterize the order structure of the set of multiple zeta-star values.  
               
The following lemma will be used in the proof of Theorem \ref{last}.
\begin{lem}\label{mgmg}
For $\alpha\in \mathbb{C},m\geq1$, we have
\[
\mathop{\int}_{[0,1]^r}[1-(1-\alpha) x_1\cdots x_r]^{m-1}dx_1\cdots dx_r=\sum_{m\geq m_1\geq \cdots\geq m_r\geq 1}\frac{1}{m m_1\cdots m_{r-1} }\alpha^{m_r-1}.
\]
\end{lem}
\noindent{\bf Proof:}
We have
\[
\begin{split}
&\;\; \mathop{\int}_{[0,1]^r}[1-(1-\alpha) x_1\cdots x_r]^{m-1}dx_1\cdots dx_r \\
&= \mathop{\int}_{[0,1]^{r-1}}\left[ \frac{1}{m}\frac{1-(1-(1-\alpha) x_2\cdots x_r)^{m}}{(1-\alpha )x_2\cdots x_r}\right]dx_2\cdots dx_r\\
&=\frac{1}{m}\sum_{m\geq m_1\geq 1}\mathop{\int}_{[0,1]^{r-1}}(1-(1-\alpha) x_2\cdots x_r)^{m_1-1}dx_2\cdots dx_r\\
&=\;\;\;\;\;  \cdots                \;\;\;\;\;\\
&=\sum_{m\geq m_1\geq \cdots\geq m_r\geq 1 }\frac{1}{m m_1\cdots m_{r-1}}\alpha^{m_r-1}.
\end{split}
\]
$\hfill\Box$\\

Now we are ready to prove Theorem \ref{last}.\\
\noindent{\bf Proof of Theorem \ref{last}:} Since 
\[
\begin{split}
&\;\;\;\;\,\frac{1}{1-x_1\cdots x_{i_1}+x_1\cdots x_{i_1}x_{i_1+1}\cdots x_{i_2}+\cdots+(-1)^rx_1\cdots x_{i_r}}\\
&= \sum_{n_1\geq 1} (x_1\cdots x_{i_1})^{n_1-1}\left[1-x_{i_1+1}\cdots x_{i_2}+\cdots+(-1)^{r-1}x_{i_1+1}\cdots x_{i_r}   \right]^{n_1-1},\\
\end{split}
\]
we have
\[
\begin{split}
&\;\;\mathop{\int}_{[0,1]^{i_r}}\frac{dx_1\cdots dx_{i_r}}{1-x_1\cdots x_{i_1}+x_1\cdots x_{i_1}x_{i_1+1}\cdots x_{i_2}+\cdots+(-1)^rx_1\cdots x_{i_r}}\\
&=\sum_{n_1\geq 1}\frac{1}{n_1^{i_1}}\mathop{\int}_{[0,1]^{i_r-i_1}}\left[1-x_{i_1+1}\cdots x_{i_2}+\cdots+(-1)^{r-1}x_{i_1+1}\cdots x_{i_r}   \right]^{n_1-1}dx_{i_1+1}\cdots dx_{i_r}\\
&=\sum_{ n_1\geq n_2\geq \cdots \geq n_{i_2-i_1+1}\geq 1  }\frac{1}{n_1^{i_1+1}n_2\cdots n_{i_2-i_1}n_{i_2-i_1+1}^{i_3-i_2}}\\
&\;\;\;\;\;\;\;\;\;\;\;\;\;\cdot\mathop{\int}_{[0,1]^{i_r-i_3}}\left[1-x_{i_3+1}\cdots x_{i_4}+\cdots+(-1)^{r-3}x_{i_3+1}\cdots x_{i_r} \right]^{n_{i_2-i_1+1}-1}.
\end{split}
\]
In the above calculation, the last identity follows from Lemma \ref{mgmg}. 
By induction, the statements $(i)$ and $(ii)$ are proved.
$\hfill\Box$\\

\begin{lem}\label{indu}
For $r\geq 1$, $k_1,\cdots, k_{r+1}\geq 2, l_1,\cdots,l_r\geq 0$,\\
$(i)$ \[
\zeta^{\star}(k_1,\{1\}^{l_1}, k_2,\{l_2\},\cdots, \{1\}^{l_r}, k_{r+1})> \zeta^{\star}(k_1,\{1\}^{l_1}, k_2,\{l_2\},\cdots, \{1\}^{l_r});\]
$(ii)$ For any $l_r>1$,
\[
\begin{split}
&\zeta^{\star}(k_1,\{1\}^{l_1}, \cdots, \{1\}^{l_{r-1}},k_r-1)>\zeta^{\star}(k_1,\{1\}^{l_1}, \cdots, \{1\}^{l_{r-1}},k_r,\{1\}^{l_r})\\
& >\zeta^{\star}(k_1,\{1\}^{l_1}, \cdots, \{1\}^{l_{r-1}},k_r,\{1\}^{l_r-1}) >\cdots> \zeta^{\star}(k_1,\{1\}^{l_1}, \cdots, \{1\}^{l_{r-1}},k_r).  \\
\end{split}
\]
\end{lem}
\noindent{\bf Proof:} $(i)$  
For $k=2r+1$, $\forall \,x_i\in(0,1), 1\leq i\leq {i_k}$, it is clear that 
\[\begin{split}
&\;\;\;\,\;1-x_1\cdots x_{i_1}+x_1\cdots x_{i_1}x_{i_1+1}\cdots x_{i_2}+\cdots+(-1)^kx_1\cdots x_{i_k}\\
&<1-x_1\cdots x_{i_1}+x_1\cdots x_{i_1}x_{i_1+1}\cdots x_{i_2}+\cdots+(-1)^{k-1}x_1\cdots x_{i_{k-1}}.\\
\end{split}\]
Thus  
\[\begin{split}
&\;\;\;\,\;\mathop{\int}_{[0,1]^{i_k}}\frac{dx_1\cdots dx_{i_k}}{1-x_1\cdots x_{i_1}+x_1\cdots x_{i_1}x_{i_1+1}\cdots x_{i_2}+\cdots+(-1)^kx_1\cdots x_{i_k}}\\
&>\mathop{\int}_{[0,1]^{i_{k-1}}}\frac{dx_1\cdots dx_{i_{k-1}}}{1-x_1\cdots x_{i_1}+x_1\cdots x_{i_1}x_{i_1+1}\cdots x_{i_2}+\cdots+(-1)^{k-1}x_1\cdots x_{i_{k-1}}}.\\ 
\end{split}\]
By Theorem \ref{last}, $(i)$, the statement $(i)$ is proved.

Similarly, the statement $(ii)$ follows from the fact:\\
For $k=2r$, $\forall \,x_i\in(0,1), 1\leq i\leq {i_k}$,
\[
\begin{split}
&\;\;\;\;1-x_1\cdots x_{i_1}+x_1\cdots x_{i_1}x_{i_1+1}\cdots x_{i_2}+\cdots+(-1)^{k-1}x_1\cdots x_{i_{k-1}}
\\
&<1-x_1\cdots x_{i_1}+x_1\cdots x_{i_1}x_{i_1+1}\cdots x_{i_2}+\cdots+(-1)^kx_1\cdots x_{i_k}\\&<\cdots \\
&\,\;\;\;\;\cdots  \;\;\;\;\;\;\;\;\;\;\;\; \cdots  \\
&<1-x_1\cdots x_{i_1}+x_1\cdots x_{i_1}x_{i_1+1}\cdots x_{i_2}+\cdots+(-1)^kx_1\cdots x_{i_{k-1}+1}.\\
\end{split}\]
$\hfill\Box$\\

Now we are ready to prove Theorem \ref{order}. By Lemma \ref{indu}, it follows that
\[
\zeta^{\star}(k_1,\cdots, k_r, k_{r+1})>\zeta^{\star}(k_1,\cdots, k_r), \;\forall\,(k_1,\cdots,k_r,k_{r+1})\in \mathcal{S},
\]
\[
\zeta^{\star}(k_1,\cdots,k_{r-1}, k_r)>\zeta^{\star}(k_1,\cdots,k_{r-1}, k_r+1, \{1\}^l),\; \forall\,(k_1,\cdots,k_r)\in \mathcal{S},l\geq 1.\]
By induction, if  \[
 (k_1,\cdots, k_{r})\succ(m_1,\cdots, m_{s}), \]
 then
 \[
  \zeta^{\star}(k_1,\cdots, k_{r})>\zeta^{\star}(m_1,\cdots, m_{s}). 
  \]
  If 
  \[
  \zeta^{\star}(k_1,\cdots, k_{r})>\zeta^{\star}(m_1,\cdots, m_{s}), 
  \]
  we must have \[
  (k_1,\cdots, k_{r})\succ(m_1,\cdots, m_{s}). \]
  Otherwise we will get a contradiction. As a result, Theorem \ref{order} is proved.
  
\section{The derived set of the set of multiple zeta-star values}
  In this section we will show that the derived set of the set of multiple zeta-star values  is the half straight line $[1,+\infty)$. Moreover, based on the order structure of the set of multiple zeta-star values, we will establish a one-to-one correspondence between the set of sequences $\mathcal{T}$ and the  half straight line $(1,+\infty)$. Lastly, we will investigate the limit of  some multiple integrals.
  
  By the results of Ohno and Wakabayashi \cite{OW}, we have 
  \[
  \zeta^{\star}(2, \{1\}^n)=(n+1)\zeta(n+2).
    \]
  
  \begin{lem}\label{min}
   For $k_1\geq 2, k_2,\cdots,k_r\geq 1$, one has  
   \[
   \zeta^{\star}(k_1,\cdots, k_{r})- \zeta^{\star}(k_1,\cdots, k_{r-1})= \sum_{n_1\geq \cdots\geq n_{r}\geq 2}   \frac{1}{n_1^{k_1}\cdots n_{r}^{k_{r}}}.  \]
     \end{lem}
  \noindent{\bf Proof:} It is clear that
 \[
 \begin{split}
 &\;\;\;\;  \zeta^{\star}(k_1,\cdots, k_{r})       \\
 &=  \sum_{n_1\geq \cdots\geq n_{r}\geq 1}   \frac{1}{n_1^{k_1}\cdots n_{r}^{k_{r}}}          \\
 &=    \sum_{\substack{n_1\geq \cdots\geq n_{r}\geq 1\\ n_r\geq 2}}   \frac{1}{n_1^{k_1}\cdots n_{r}^{k_{r}}}  +  \sum_{\substack{n_1\geq \cdots\geq n_{r}\geq 1\\ n_r=1}}   \frac{1}{n_1^{k_1}\cdots n_{r}^{k_{r}}}      \\
 &=  \sum_{n_1\geq \cdots\geq n_{r}\geq 2}   \frac{1}{n_1^{k_1}\cdots n_{r}^{k_{r}}}+   \zeta^{\star}(k_1,\cdots, k_{r-1}).     \\
 \end{split}
 \]
  $\hfill\Box$\\

 \begin{lem}\label{plus}
 For ${\bf k}=(k_1,\cdots,k_r,\cdots)\in\mathcal{T}$, define
 \[
 \Delta_r=\zeta^{\star}(k_1,\cdots,k_r,1)-\zeta^{\star}(k_1,\cdots,k_r).
  \]
  Then 
  \[
 \mathop{\mathrm{lim}}_{r\rightarrow +\infty}  \Delta_r=0. \]
  \end{lem}
  \noindent{\bf Proof:}
  Since ${\bf k}=(k_1,\cdots,k_r,\cdots)\in\mathcal{T}$, there is an $s\geq 2$  such that $k_s\geq 2$.
  By Theorem \ref{order}, for $m\geq 1$,
  \[
 \zeta^{\star}(k_1,\cdots,k_{s-1},k_s-1) >\zeta^{\star}(k_1,\cdots,k_s,{\{1\}}^{m+1})>\zeta^{\star}(k_1,\cdots,k_s,{\{1\}}^m) >\zeta^{\star}(k_1,\cdots,k_s).  \]
 So the limit 
  \[
  \mathop{\mathrm{lim}}_{m\rightarrow +\infty} \zeta^{\star}(k_1,\cdots, k_{s},\{1\}^m) \]
  exists and 
  \[
  \mathop{\mathrm{lim}}_{m\rightarrow +\infty}\left( \zeta^{\star}(k_1,\cdots, k_{s},\{1\}^{m}) - \zeta^{\star}(k_1,\cdots, k_{s},\{1\}^{m-1}) \right)=0. \]
  By Lemma \ref{min}, it follows that
  \[
   \mathop{\mathrm{lim}}_{m\rightarrow +\infty}\sum_{n_1\geq \cdots\geq n_{s}\geq \cdots \geq n_{s+m}\geq 2}   \frac{1}{n_1^{k_1}\cdots n_{s}^{k_{s}}n_{s+1}\cdots n_{s+m}}=0.     \]
   For $r>s$, we have 
   \[
   0<\Delta_r=\sum_{n_1\geq \cdots\geq n_{r}\geq n_{r+1}\geq 2}   \frac{1}{n_1^{k_1}\cdots n_{r}^{k_{r}}n_{r+1}}\leq \sum_{n_1\geq \cdots\geq n_{r}\geq n_{r+1}\geq 2}   \frac{1}{n_1^{k_1}\cdots n_{s}^{k_{s}}n_{s+1}\cdots n_{r+1}}. \]
   Thus 
   \[
 \mathop{\mathrm{lim}}_{r\rightarrow +\infty}  \Delta_r=0. \]
  $\hfill\Box$\\

  \begin{rem}
  In contrast to Lemma \ref{plus}, the limit 
  \[
    \mathop{\mathrm{lim}}_{m\rightarrow +\infty}\left( \zeta^{\star}(2,\{1\}^{m}) - \zeta^{\star}(2,\{1\}^{m-1}) \right) \] is not zero. In fact, we have 
    \[
    \begin{split}
    &\;\;\;\;   \mathop{\mathrm{lim}}_{m\rightarrow +\infty}\left( \zeta^{\star}(2,\{1\}^{m}) - \zeta^{\star}(2,\{1\}^{m-1}) \right)             \\
    &=  \mathop{\mathrm{lim}}_{m\rightarrow +\infty} \left( (m+1)\zeta(m+2)-m\zeta(m+1)    \right)                    \\
    &=1.
    \end{split}
    \]
  \end{rem}

   \begin{lem}\label{sep}
 $(i)$ 
  \[
 \mathop{\mathrm{lim}}_{n\rightarrow +\infty} \zeta^{\star}(2, \{1\}^n)=+\infty.
 \]
  $(ii)$ For $k_1\geq 2, k_2,\cdots,k_r\geq 1$, one has
 \[
  \mathop{\mathrm{lim}}_{n\rightarrow +\infty} \zeta^{\star}(k_1,\cdots, k_{r-1},k_r+1,\{1\}^n)= \zeta^{\star}(k_1,\cdots, k_{r-1},k_r), \]
   \[
  \mathop{\mathrm{lim}}_{n\rightarrow +\infty} \zeta^{\star}(k_1,\cdots, k_{r-1},k_r,n)=\zeta^{\star}(k_1,\cdots, k_{r-1},k_r). \]
 \end{lem} 
\noindent{\bf Proof:}
$(i)$  The statement $(i)$ follows immediately from the formula 
\[
\zeta^{\star}(2, \{1\}^n)=(n+1)\zeta(n+2).\]
Here we give a direct proof based on Theorem \ref{last}.
By Theorem \ref{last}, it follows that
\[
\begin{split}
&\;\;\;\;\zeta^{\star}(2, \{1\}^n)\\
&=\mathop{\int}_{[0,1]^{n+2}}\frac{dx_1dx_2\cdots dx_{n+2}}{1-x_1+x_1x_2\cdots x_{n+2}}\\
&= \mathop{\int}_{[0,1]^{n+2}}    \sum_{m\geq 1}x_1^{m-1}(1-x_2\cdots x_{n+2})^{m-1}      dx_1dx_2\cdots dx_{n+2}  \\
&=  \sum_{m\geq 1}\frac{1}{m}\mathop{\int}_{[0,1]^{n+1}}  (1-x_2\cdots x_{n+2})^{m-1}      dx_2\cdots dx_{n+2} .\end{split}
\]
For any $m\geq 1$, one can check that
\[
\begin{split}
&\;\;\; \mathop{\mathrm{lim}}_{n\rightarrow +\infty} \mathop{\int}_{[0,1]^{n+1}}  (1-x_2\cdots x_{n+2})^{m-1}      dx_2\cdots dx_{n+2}\\
&= \mathop{\mathrm{lim}}_{n\rightarrow +\infty} \mathop{\int}_{[0,1]^{n+1}}  \sum_{1\leq p\leq m}(-1)^{p-1} \binom{m-1}{p-1}(x_2\cdots x_{n+2})^{p-1}     dx_2\cdots dx_{n+2}  \\
&=  \mathop{\mathrm{lim}}_{n\rightarrow +\infty} \sum_{1\leq p\leq m}\binom{m-1}{p-1}\frac{(-1)^{p-1}}{p^{n+1}}\\
&=1.
\end{split}
 \]
 As a result, for any $\varepsilon>0, M\geq 1$, there is an $N$ such that for $n>N$, 
 \[
 \zeta^{\star}(2, \{1\}^n)\geq \sum_{M\geq m\geq 1}\frac{1}{m}-\varepsilon.
  \]
  Since 
  \[
    \mathop{\mathrm{lim}}_{M\rightarrow +\infty} \sum_{M\geq m\geq 1}\frac{1}{m}=+\infty,      \]
    we have  \[
 \mathop{\mathrm{lim}}_{n\rightarrow +\infty} \zeta^{\star}(2, \{1\}^n)=+\infty.
 \]
 
 $(ii)$ By Theorem \ref{order}, $ \forall\, n\geq 1$,
  \[
 \zeta^{\star}(k_1,\cdots, k_{r-1},k_r+1,\{1\}^{n+1})>\zeta^{\star}(k_1,\cdots, k_{r-1},k_r+1,\{1\}^n),  \]
 \[\zeta^{\star}(k_1,\cdots, k_{r-1},k_r+1,\{1\}^n)<\zeta^{\star}(k_1,\cdots, k_{r-1},k_r). \]
Thus the limit  $$\mathop{\mathrm{lim}}_{n\rightarrow +\infty} \zeta^{\star}(k_1,\cdots, k_{r-1},k_r+1,\{1\}^n)$$ exists.
By Lemma \ref{min}, we have 
\[
\begin{split}
&\;\;\;\;  \mathop{\mathrm{lim}}_{n\rightarrow +\infty} \zeta^{\star}(k_1,\cdots, k_{r-1},k_r+1,\{1\}^n)    \\
&= 1+\sum_{n_2\geq 2}\frac{1}{n_1^{k_1}}+\cdots+ \sum_{n_1\geq \cdots\geq n_{r-1}\geq 2}\frac{1}{n_1^{k_1}\cdots n_{r-1}^{k_{r-1}}}           
+ \sum_{n_1\geq \cdots\geq n_r\geq 2}\frac{1}{n_1^{k_1}\cdots n_{r-1}^{k_{r-1}}n_r^{k_r+1}}\Bigg{(} 1\\
&+\sum_{n_r\geq n_{r+1}\geq 2}\frac{1}{n_{r+1}}+\cdots+  \sum_{n_r\geq n_{r+1}\geq \cdots\geq n_{r+s}\geq 2}\frac{1}{n_{r+1}\cdots n_{r+s}} +\cdots   \Bigg{)}            \\
&=  1+\sum_{n_2\geq 2}\frac{1}{n_1^{k_1}}+\cdots+ \sum_{n_1\geq \cdots\geq n_{r-1}\geq 2}\frac{1}{n_1^{k_1}\cdots n_{r-1}^{k_{r-1}}} \\
&+   \sum_{n_1\geq \cdots\geq n_r\geq 2}\frac{1}{n_1^{k_1}\cdots n_{r-1}^{k_{r-1}}n_r^{k_r+1}} \prod_{n_r\geq l\geq 2}\left(1+\frac{1}{l}+\cdots+\frac{1}{l^j}+\cdots\right)         \\
\end{split}
\]

  Since 
  \[
   \prod_{n_r\geq l\geq 2}\left(1+\frac{1}{l}+\cdots+\frac{1}{l^j}+\cdots\right)  =\prod_{n_r\geq l\geq 2}\frac{l}{l-1}=n_r,
     \]
     it follows that
       \[
\begin{split}
&\;\;\;\;  \mathop{\mathrm{lim}}_{n\rightarrow +\infty} \zeta^{\star}(k_1,\cdots, k_{r-1},k_r+1,\{1\}^n)    \\
&= 1+\sum_{n_2\geq 2}\frac{1}{n_1^{k_1}}+\cdots+ \sum_{n_1\geq \cdots\geq n_{r-1}\geq 2}\frac{1}{n_1^{k_1}\cdots n_{r-1}^{k_{r-1}}}           
+ \sum_{n_1\geq \cdots\geq n_r\geq 2}\frac{1}{n_1^{k_1}\cdots n_{r-1}^{k_{r-1}}n_r^{k_r}}\\
&=\zeta^{\star}(k_1,\cdots, k_{r-1},k_r).
\end{split}
\]
Similarly, we have 
\[
\begin{split}
&\;\;\;\;  \mathop{\mathrm{lim}}_{n\rightarrow +\infty} \zeta^{\star}(k_1,\cdots, k_{r-1},k_r,n)    \\
&= \zeta^{\star}(k_1,\cdots, k_{r-1},k_r)+ \mathop{\mathrm{lim}}_{n\rightarrow +\infty} \sum_{n_1\geq n_r\geq n_{r+1}\geq 2}  \frac{1}{n_1^{k_1}\cdots n_r^{k_r} n_{r+1}^n }  \\
&= \zeta^{\star}(k_1,\cdots, k_{r-1},k_r)+  \sum_{n_1\geq n_r\geq n_{r+1}\geq 2} \mathop{\mathrm{lim}}_{n\rightarrow +\infty} \frac{1}{n_1^{k_1}\cdots n_r^{k_r} n_{r+1}^n }  \\
&= \zeta^{\star}(k_1,\cdots, k_{r-1},k_r).
\end{split}
\] 
 $\hfill\Box$\\    

 \begin{lem}\label{mz} For ${\bf k}=(k_1,\cdots,k_r,\cdots)\in\mathcal{T}$ and $(k_1,\cdots, k_t)\neq (2,\{1\}^{t-1})$,  define
  \[
  \Delta_{k_1,\cdots, k_r}= \sum_{n_1\geq \cdots \geq n_r\geq 1}\frac{n_r-1}{n_1^{k_1}n_2^{k_2}\cdots n_r^{k_r}    }
  \]  
  for $r\geq t$.  Then  \[
 \mathop{\mathrm{lim}}_{r\rightarrow +\infty}  \Delta_{k_1,\cdots,k_r}=0. \]
   \end{lem}
    \noindent{\bf Proof:}
    For $r\geq t$, 
    it is clear that
    \[
    \begin{split}
    &\;\;\;\;   \Delta_{k_1,\cdots, k_r}\\
    &\leq   \sum_{n_1\geq \cdots n_r\geq 1}\frac{1}{n_1^{k_1}\cdots n_t^{k_t}} \cdot \frac{n_r-1}{n_{t+1}\cdots n_r      }    .                                  \\
        \end{split}
    \]
    If $k_i\geq 2, k_{i+1}=\cdots =k_t=1$ for some $i\geq r$, then
    
       \[
    \begin{split}
    &\;\;\;\;   \Delta_{k_1,\cdots, k_r}\\
    &\leq   \sum_{n_1\geq \cdots n_r\geq 1}\frac{1}{n_1^{k_1}\cdots n_i^{k_i}} \cdot \frac{n_r-1}{n_{i+1}\cdots n_r      }    \\                                \\
    &\leq \sum_{n_1\geq \cdots n_i\geq 1}\frac{1}{n_1^{k_1}\cdots n_i^{k_i}} \cdot  \left( \sum_{n_i\geq n_{i+1}\geq \cdots \geq n_r\geq 1}    \frac{n_r-1}{n_{i+1}\cdots n_r      }    \right)  \\
    &\leq \sum_{n_1\geq \cdots n_i\geq 1}\frac{1}{n_1^{k_1}\cdots n_i^{k_i}} \cdot  \left[ \sum_{n_i\geq n_{i+1}\geq \cdots \geq n_r\geq 1} \left(\frac{1}{n_{i+1}\cdots n_{r-1}      }  - \frac{1}{n_{i+1}\cdots n_r      }  \right)  \right]  \\
       &\leq \sum_{n_1\geq \cdots n_i\geq 1}\frac{1}{n_1^{k_1}\cdots n_i^{k_i}} \cdot  \left(n_i- \sum_{n_i\geq n_{i+1}\geq \cdots \geq n_r\geq 1} \frac{1}{n_{i+1}\cdots n_r      }  \right)  \\
       &\leq \zeta^\star(k_1,\cdots, k_{i-1}, k_i-1)-\zeta^\star(k_1,\cdots, k_{i-1},k_i,\{1\}^{r-i}).\\
            \end{split}
    \]
    Here the fourth inequality follows from the indentity
    \[
    n_i=\sum_{n_i\geq n_{i+1}\geq \cdots \geq n_r\geq 1}\frac{1}{n_{i+1}\cdots n_{r-1}      }      \]
    for every $n_i\geq 1$.
    By Lemma \ref{sep}, $(ii)$, we have 
    \[ \mathop{\mathrm{lim}}_{r\rightarrow +\infty}  \Delta_{k_1,\cdots,k_r}=0. \]
     $\hfill\Box$\\  

\noindent{\bf Proof of Theorem \ref{rep}:}
$(i)$ For any $(k_1,k_2,\cdots,k_r,\cdots)\in \mathcal{T}$, denote by $$s= \mathop{min}_{j}\big{\{}j\,|\, j\geq 2\; \mathrm{and}\; k_j\geq 2\big{\}}.$$
By Theorem \ref{order}, it follows that $\forall\,r\geq 1$,
\[
\zeta^\star(k_1,k_2,\cdots,k_r,k_{r+1})>\zeta^\star(k_1,k_2,\cdots,k_r), 
\]
\[
\zeta^\star(k_1,k_2,\cdots,k_r)<\zeta^\star(k_1,k_2,\cdots,k_{s-1},k_s-1).\]
As a result, the limit  $$\mathop{\mathrm{lim}}_{r\rightarrow +\infty} \zeta^{\star}(k_1,\cdots, k_{r-1},k_r)$$ exists and 
$$\mathop{\mathrm{lim}}_{r\rightarrow +\infty} \zeta^{\star}(k_1,\cdots, k_{r-1},k_r)>1.$$

$(ii)$ For ${\bf p}=(p_1,p_2,\cdots,p_r,\cdots,), {\bf q}=(q_1,q_2,\cdots, q_r,\cdots)\in \mathcal{T}$ and 
\[\eta(\bf p)=\eta(\bf q),\]
we want to show that ${\bf p}={\bf q}$. If ${\bf p}\neq {\bf q}$, by exchanging the values of $\bf p$ and $\bf q$ if necessary, there is an $s\geq 1$ such that
\[
p_i=q_i, 1\leq i\leq s-1,\]
\[
 p_s<q_s.
\]
By Theorem \ref{order}, we have
\[
\zeta^{\star}(p_1,\cdots,p_{s-1},p_s)\geq \zeta^{\star}(q_1,\cdots,q_{s-1},q_s-1)>\zeta^{\star}(q_1,\cdots,q_{s-1},q_s)\]
So $\forall\,m\geq 1$,
\[
\mathop{\mathrm{lim}}_{r\rightarrow +\infty} \zeta^{\star}(p_1,\cdots, p_{r-1},p_r)> \zeta^{\star}(q_1,\cdots,q_{s-1},q_s-1)>\zeta^{\star}(q_1,\cdots,q_{s-1},q_s,\{1\}^m). \]
Since 
\[
\zeta^{\star}(q_1,\cdots,q_{s-1},q_s,\{1\}^m)\geq \zeta^{\star}(q_1,\cdots,q_{s-1},q_s,q_{s+1},\cdots,q_{s+m}),
\]
it follows that
\[
 \zeta^{\star}(q_1,\cdots,q_{s-1},q_s-1)\geq \mathop{\mathrm{lim}}_{r\rightarrow +\infty} \zeta^{\star}(q_1,\cdots, q_{r-1},q_r) .\]
 As a result, 
 \[
 \mathop{\mathrm{lim}}_{r\rightarrow +\infty} \zeta^{\star}(p_1,\cdots, p_{r-1},p_r)>\mathop{\mathrm{lim}}_{r\rightarrow +\infty} \zeta^{\star}(q_1,\cdots, q_{r-1},q_r)  .
 \] 
 There is a contradiction. So the map $\eta$ is injective. 
 
 For $x\in (1,+\infty)$, if $x\in \mathcal{Z}^\star$, then
 \[
 x=\zeta^\star(k_1,k_2,\cdots, k_r)
 \]
 for some $(k_1,k_2,\cdots,k_r)\in\mathcal{S}$. By Lemma \ref{sep}, $(ii)$, 
 \[
 x= \mathop{\mathrm{lim}}_{t\rightarrow +\infty} \zeta^{\star} (p_1,p_2,\cdots, p_t),
   \]
   where $$(p_1,p_2,\cdots, p_t,\cdots)=  (k_1,\cdots, k_{r-1},k_r+1,1,\cdots,1,\cdots).$$
   For $k_1\geq 2$, define 
   \[
   {Z}_{k_1}=\begin{cases} \big{(} \zeta(k_1),\zeta(k_1-1)     \big{ ]}, &  k_1\geq 3; \\
   \big{(}\zeta(2), +\infty\big{)}, & k_1=2. 
   \end{cases}
   \]
   For $r\geq 2$, $k_1\geq 2, k_2,\cdots, k_r\geq 1$, define 
 \[
 Z_{k_1,\cdots, k_r}=\begin{cases} \big{(} \zeta^\star(k_1,\cdots, k_{r-1}, k_r),\zeta^\star(k_1,\cdots, k_{r-1}, k_r-1)     \big{ ]}, &  k_r\geq 2; \\
   \big{(}\zeta^\star(k_1,\cdots, k_{r-1}, k_r ), \zeta^\star(k_1,\cdots, k_{i-1}, k_i-1 ) \big{]}, & k_i\geq 2, k_{i+1}=\cdots =k_r=1. 
   \end{cases}
   \]
   Here we assume that $\zeta(1)=+\infty$.
   By Theorem \ref{order} and  Lemma \ref{sep}, it is easy to check that
   \[
 {Z}_{p_1,\cdots, p_r}\cap Z_{q_1,\cdots, q_r} =\emptyset , \mathrm{for}\,(p_1,\cdots,p_r)\neq (q_1,\cdots, q_r),r\geq 1,\tag{1}     \]
   \[
   (1,+\infty)= \bigcup_{\substack{  k_1\geq 2,k_2,\cdots,k_r\geq 1            \\                  }              } Z_{k_1,\cdots,k_r},  \tag{2} \]
   \[
Z_{k_1,\cdots, k_t}= \bigcup_{k_{t+1},\cdots, k_{t+s}\geq 1} Z_{k_1,\cdots, k_t, k_{t+1},\cdots, k_{t+s}} . \tag{3}  \]

   If  $x\notin \mathcal{Z}^\star$, by  the formulas $(1), (2)$ and $(3)$ and Lemma \ref{sep}, $(i)$, there is a 
   $${\bf k}=(k_1,k_2, \cdots, k_r,\cdots)\in  \mathcal{T}$$  such that
   \[
   x\in Z_{k_1,k_2,\cdots,k_r} \tag{4}
   \]
   for all  $r\geq 1$.
   If $(k_1,\cdots, k_r)\neq (2,\{1\}^{r-1})$ for $r\geq t$, then 
   $$x<\zeta^\star(2,\{1\}^{t-1})$$
   and 
    \[
    \begin{split}
    &\;\;\;\;m\left( Z_{k_1,\cdots,k_r}      \right)\\
    &=\zeta^\star(k_1,\cdots, k_{r-1}, k_r-1)-\zeta^{\star}(k_1,\cdots, k_{r-1},k_r)                                              \\
    &= \sum_{n_1\geq \cdots \geq n_{r-1}\geq n_r\geq 1}\frac{n_r-1}{ n_1^{k_1}\cdots n_{r-1}^{k_{r-1}} n_r^{k_r}         }                             \\
      &=     \sum_{n_1\geq \cdots \geq n_{r-1}\geq n_r\geq 2}\frac{n_r-1}{ n_1^{k_1}\cdots n_{r-1}^{k_{r-1}} n_r^{k_r}         }                    \\
 \end{split}
  \]
  for $k_r\geq 2$
  and 
   \[
    \begin{split}
    &\;\;\;\;m\left( Z_{k_1,\cdots,k_r}      \right)\\
    &=\zeta^\star(k_1,\cdots, k_{i-1}, k_i-1)-\zeta^{\star}(k_1,\cdots, k_{r-1},k_r)                                              \\
    &=\sum_{n_1\geq \cdots \geq n_{i-1}\geq n_i\geq 1}\frac{1}{n_1^{k_1} \cdots n_{i-1}^{k_{i-1}} n_i^{k_i-1}   }-  \sum_{n_1\geq \cdots \geq n_{r-1}\geq n_r\geq 1}\frac{1}{n_1^{k_1} \cdots n_{r-1}^{k_{r-1}} n_r^{k_r}   }         \\        
    &=   \sum_{n_1\geq \cdots \geq n_{i-1}\geq n_i\geq 1}\frac{1}{n_1^{k_1} \cdots n_{i-1}^{k_{i-1}} n_i^{k_i}   } \left(  \sum_{n_i\geq n_{i+1}\geq \cdots \geq n_r\geq 1} \frac{1}{n_{i+1}\cdots n_{r-1}}\right)\\
    &\;\;\;\;  -  \sum_{n_1\geq \cdots \geq n_{r-1}\geq n_r\geq 1}\frac{1}{n_1^{k_1} \cdots n_{r-1}^{k_{r-1}} n_r^{k_r}   }         \\        
    &= \sum_{n_1\geq \cdots \geq n_{r-1}\geq n_r\geq 1}\frac{n_r-1}{ n_1^{k_1}\cdots n_{r-1}^{k_{r-1}} n_r^{k_r}         }       \\
      &=     \sum_{n_1\geq \cdots \geq n_{r-1}\geq n_r\geq 2}\frac{n_r-1}{ n_1^{k_1}\cdots n_{r-1}^{k_{r-1}} n_r^{k_r}         }         \\
   \end{split}
  \]
     for $k_i\geq 2, k_{i+1}=\cdots =k_r=1$.
     Here the third equality follows from the fact that
     \[
     n_i=   \sum_{n_i\geq n_{i+1}\geq \cdots \geq n_r\geq 1} \frac{1}{n_{i+1}\cdots n_{r-1}}     \]
     for every $n_i\geq 1$.
   From $(4)$ one has 
   \[
    \zeta^\star(k_1,\cdots, k_{r-1}, k_r)<x\leq    \zeta^\star(k_1,\cdots, k_{r-1}, k_r)+  m\left( Z_{k_1,\cdots,k_r}      \right)   \]
    for $r\geq t+1$.
    By Lemma \ref{mz}, it follows that
        $$x=\mathop{\mathrm{lim}}_{r\rightarrow +\infty} \zeta^{\star}(k_1,\cdots, k_{r-1},k_r).$$ 
        The fact that $\eta$ preserves the order structure follows immediately from the above construction.  
        $(iii)$ The statement $(iii)$ follows immediately from the statement $(ii)$.  
           $\hfill\Box$\\           
        
        For ${\bf i}=(i_1,i_2,\cdots, i_r,\cdots)\in\mathcal{A}$, by Theorem \ref{last}, there is a ${\bf k}=(k_1,\cdots,k_r,\cdots)\in\mathcal{T}$ such that 
        \[
\begin{split}
&\;\;\mathop{\int}_{[0,1]^{i_r}}\frac{dx_1\cdots dx_{i_r}}{1-x_1\cdots x_{i_1}+x_1\cdots x_{i_2}+\cdots+(-1)^rx_1\cdots x_{i_r}}\\
&= \begin{cases}      
\zeta^{\star}(k_1,\cdots, k_{l_r-1},k_{l_r}-1),& \mathrm{for\;some}\; l_r>\frac{r}{2}, k_{l_r}\geq 2,r\,\mathrm{odd} ;  \\
 \zeta^{\star}(k_1,\cdots, k_{l_r-1},k_{l_r}),& \mathrm{for\;some}\; l_r>\frac{r}{2}, r\,\mathrm{even}.  \\
\end{cases}
\end{split}
\]
In fact, on can check that 
\[\tiny
\begin{split}
&\;\;\;\;(k_1, \cdots,k_r,\cdots)\\
&=(i_{1}+1,\{1\}^{i_2-i_1-1}, i_3-i_2+1,\{1\}^{i_4-i_3-1},\cdots, \cdots, i_{{2k-1}}- i_{{2k-2}} +1, \{1\}^{i_{{2k}}-i_{{2k-1}}-1 },\cdots      ).\\ 
\end{split} \tag{5}\]
Since 
\[
  \zeta^{\star}(k_1,\cdots, k_{l_r-1},1)\geq \zeta^{\star}(k_1,\cdots, k_{l_r-1},k_{l_r}-1)> \zeta^{\star}(k_1,\cdots, k_{l_r-1},k_{l_r}),
\]
by Lemma \ref{min}, 
 \[
 \begin{split}
&\;\;\;\;\mathop{\mathrm{lim}}_{r\rightarrow +\infty} \mathop{\int}_{[0,1]^{i_r}}\frac{dx_1dx_2\cdots dx_{i_r}}{1-x_1x_2\cdots x_{i_1}+ x_1x_2\cdots x_{i_2}  +\cdots+(-1)^rx_1x_2\cdots x_{i_r}    }\\
&=   \mathop{\mathrm{lim}}_{r\rightarrow +\infty}    \zeta^{\star}(k_1,\cdots, k_{l_r-1},k_{l_r})     \\
&= \mathop{\mathrm{lim}}_{r\rightarrow +\infty}    \zeta^{\star}(k_1,\cdots, k_{r-1},k_{r}) .
\end{split}
 \]        
 Thus Corollary \ref{alt}, $(i)$ is proved.
 
 The map $\xi$ is injective follows from
the formula $(5)$ and $\eta$ is injective.
By the formula $(5)$, for 
\[
x=\mathop{\mathrm{lim}}_{r\rightarrow +\infty}    \zeta^{\star}(k_1,\cdots, k_{r-1},k_{r})\in (1,+\infty),
\]
$x\in \xi(\mathcal{A})$ if and only if 
\[
\sharp \{i\,|\, k_i\geq 2\}=+\infty. \tag{6}
\]
By Theorem \ref{rep} and  Lemma \ref{sep}, $(ii)$, the condition $(6)$ is equivalent to 
\[
x\in (1,+\infty)-\mathcal{Z}^{\star}.
\]
In a word, $\xi(\mathcal{A})=(1,+\infty)-\mathcal{Z}^{\star}$. Corollary \ref{alt}, $(ii)$ is proved.
\begin{rem}
 Define 
 \[
 \mathcal{T}_{\mathbb{R}}=\Big{\{} (x_1,x_2, \cdots, x_r,\cdots)\,\Big{|}\,x_1\geq 2,x_i\geq 1, \forall \,i\geq 2, x_s\geq 2 \;for \; some \; s\geq 2  , x_i \in \mathbb{R}    \Big{\}}.
 \]
 By exactly the same analysis, one can show that 
 for any $${\bf x}=(x_1,x_2, \cdots, x_r,\cdots)\in  \mathcal{T}_{\mathbb{R}},$$ the limit 
  \[
 I=\mathop{\mathrm{lim}}_{r\rightarrow +\infty}\zeta^{\star}(x_1,x_2,\cdots,x_r).
 \]
 exists and  $I>1$. 
 Furthermore, the map 
 \[
 \eta_{\mathbb{R}}: \mathcal{T}_{\mathbb{R}}\rightarrow (1, +\infty),
  \]
  \[
  {\bf x}=(x_1,x_2, \cdots, x_r,\cdots) \mapsto \mathop{\mathrm{lim}}_{r\rightarrow +\infty}\zeta^{\star}(x_1,x_2,\cdots,x_r)  \]
 is surjective.  For $y\in (1,+\infty)$, its inverse image $\eta_{\mathbb{R}}^{-1}(y)$ remains mysterious.
 \end{rem}

\section{The distribution of the integer sequence }

For $x\in (1,+\infty)$, there is a ${\bf k}=(k_1,\cdots,k_r,\cdots)\in\mathcal{T}$ such that $\eta({\bf k})=x$. In this section, for some special subsets of $\mathcal{T}$, we investigate their images under the map $\eta$. Furthermore, we give the Hausdorff dimensions for some of the images. For the definition of  the Hausdorff dimension of a subset of $\mathbb{R}^n$, the reference is \cite{fal}.

  \begin{lem}\label{leb} For $s\geq 1$, 
\[  \sum_{n_1\geq \cdots\geq n_{s}\geq n_{s+1}\geq 2} \frac{1}{n_1(n_1-1)\cdots n_{s}(n_{s}-1) n_{s+1} }\leq \frac{1}{2^{s+1}} \prod_{l\geq 3}\frac{1}{1-\frac{2}{l(l-1)}}+\zeta^\star(\{2\}^s,1)-\zeta^\star(\{2\}^s).   \]  \end{lem}
\noindent{\bf Proof:} 
\[
\begin{split}
&\;\;\;\; \sum_{n_1\geq \cdots\geq n_{s}\geq n_{s+1}\geq 2} \frac{1}{n_1(n_1-1)\cdots n_{s}(n_{s}-1) n_{s+1} }    \\
&=\frac{1}{2^{s+1}}+ \frac{1}{2^{s+1-j}}\sum_{1\leq j\leq s}\sum_{n_1\geq \cdots\geq n_j\geq 3} \frac{1}{n_1(n_1-1)\cdots n_j(n_j-1)}  \\
&\;\;\;\; +\sum_{n_1\geq \cdots\geq n_{s+1}\geq 3}  \frac{1}{n_1(n_1-1)\cdots n_{s}(n_{s}-1) n_{s+1} }     \\
&<\frac{1}{2^{s+1}}\prod_{l\geq 3} \left[  1+\frac{2}{l(l-1)}+\cdots+\left( \frac{2}{l(l-1)}  \right)^j+  \cdots   \right]+
\sum_{n_1\geq \cdots\geq n_{s+1}\geq 2}\frac{1}{n_1^2\cdots n_s^2 n_{s+1}}\\
&\leq \frac{1}{2^{s+1}} \prod_{l\geq 3}\frac{1}{1-\frac{2}{l(l-1)}}+\zeta^\star(\{2\}^s,1)-\zeta^\star(\{2\}^s).\\
\end{split}
\]
 $\hfill\Box$\\

   \begin{lem}\label{dbou}
 For $p,s\geq 2$ and $(k_1,\cdots, k_t)\neq (2,\{1\}^{t-1})$, one has 
 \[
 \mathop{\mathrm{lim}}_{s\rightarrow+\infty} \sum_{n_1\geq \cdots \geq n_{t+s}\geq 2} \frac{1}{n_1^{k_1}\cdots n_t^{k_t}} \cdot \frac{1}{(n_{t+1}-1)\cdots (n_{t+s-1}-1)} \cdot \left(1-\frac{1}{n_{t+1}^p}\right)\cdots   \left(1-\frac{1}{n_{t+s}^p}\right) =0.    \]
 \end{lem} 
  \noindent{\bf Proof:}
  For $s\geq 2$, define 
  \[
  A_s=  \sum_{n_1\geq \cdots \geq n_{t+s}\geq 2} \frac{1}{n_1^{k_1}\cdots n_t^{k_t}} \cdot \frac{1}{(n_{t+1}-1)\cdots (n_{t+s-1}-1)} \cdot \left(1-\frac{1}{n_{t+1}^p}\right)\cdots   \left(1-\frac{1}{n_{t+s}^p}\right)        ,
  \]
  \[
  B_s=    \sum_{n_1\geq \cdots \geq n_{t+s-1}\geq 2} \frac{1}{n_1^{k_1}\cdots n_t^{k_t}} \cdot \frac{1}{(n_{t+1}-1)\cdots (n_{t+s-1}-1)} \cdot \left(1-\frac{1}{n_{t+1}^p}\right)\cdots   \left(1-\frac{1}{n_{t+s-1}^p}\right)        ,
  \]
  \[
  C_s=   \sum_{n_1\geq \cdots \geq n_{t+s}\geq 3} \frac{1}{n_1^{k_1}\cdots n_t^{k_t}} \cdot \frac{1}{(n_{t+1}-1)\cdots (n_{t+s-1}-1)} \cdot \left(1-\frac{1}{n_{t+1}^p}\right)\cdots   \left(1-\frac{1}{n_{t+s}^p}\right)              .
  \]
  One has 
  \[
  \begin{split}
  &\;\;\;\; A_s\\
  &= \left(  \sum_{\substack{n_1\geq \cdots \geq n_{t+s}\geq 2\\ n_{t+s}=2        }}  + \sum_{\substack{n_1\geq \cdots \geq n_{t+s}\geq 2\\ n_{t+s}\geq 3        }}    \right)  \\
  &\;\;\;\;    \frac{1}{n_1^{k_1}\cdots n_t^{k_t}} \cdot \frac{1}{(n_{t+1}-1)\cdots (n_{t+s-1}-1)} \cdot \left(1-\frac{1}{n_{t+1}^p}\right)\cdots   \left(1-\frac{1}{n_{t+s}^p}\right)         \\
  &= \left(1-\frac{1}{2^p}\right)    B_s+C_s          \\
  \end{split}
  \]
  Since 
  \[
  \begin{split}
  &\;\;\;\; C_s\\
  &<   \sum_{n_1\geq \cdots \geq n_{t+s}\geq 3} \frac{1}{n_1^{k_1}\cdots n_t^{k_t}} \cdot \frac{1}{(n_{t+1}-1)\cdots (n_{t+s-1}-1)}          \\
  &< \sum_{n_1\geq \cdots \geq n_{t+s}\geq 3} \frac{1}{(n_1-1)^{k_1}\cdots (n_t-1)^{k_t}} \cdot \frac{1}{(n_{t+1}-1)\cdots (n_{t+s-1}-1)}          \\
  &\leq  \sum_{n_1\geq \cdots \geq n_{t+s}\geq 2} \frac{1}{n_1^{k_1}\cdots n_t^{k_t}} \cdot \frac{1}{n_{t+1}\cdots n_{t+s-1}}  \\
  &\leq \sum_{n_1\geq \cdots\geq n_{t+s-1}\geq 2}\frac{1}{n_1^{k_1}\cdots n_t^{k_t}} \cdot \frac{n_{t+s-1}-1}{n_{t+1}\cdots n_{t+s-1}}, \\
 \end{split}
  \]
  By Lemma \ref{mz}, one has 
  \[
  \mathop{\mathrm{lim}}_{s\rightarrow +\infty} C_s=0.
  \]
  By the same trick, one can show that 
  \[
  B_s=\left(1-\frac{1}{2^p}\right)    B_{s-1}+D_s      \]
  for some $D_s>0$ 
  and 
    \[
  \mathop{\mathrm{lim}}_{s\rightarrow +\infty} D_s=0.
  \]
  Thus 
  \[
    \mathop{\mathrm{lim}}_{s\rightarrow +\infty} B_s= \left(1-\frac{1}{2^p}\right)^i    \mathop{\mathrm{lim}}_{s\rightarrow +\infty}   B_{s-i}=0.  \]
    As a result,
    \[
      \mathop{\mathrm{lim}}_{s\rightarrow +\infty} A_s=0.    \]
     $\hfill\Box$\\

 \noindent{\bf Proof of Theorem \ref{zero}:}           $(i)$ 
 Since $$\mathcal{T}_{p+1}\subseteq \mathcal{T}_{p}, p\geq 2,$$  it suffices to show that 
         \[
 \mu(\eta( \mathcal{T}_2))=0. \]     
 For $s\geq 1$, define 
 \[
 M_s=\{(k_1,\cdots,k_r,\cdots)\in\mathcal{T}\,|\, k_1,k_2,\cdots,k_{s+1}\geq 2.    \}
 \]
 It is clear that
 \[
 \mathcal{T}_2\subseteq M_s, \forall\,  s\geq 1.
 \]
 For $x\in(1,+\infty)$, by construction of the map $\eta$, $x\in \eta(M_s)$ if and only if 
 \[
 \zeta^{\star}(k_1,\cdots,k_{s+1})<x<\zeta^{\star}(k_1,\cdots,k_{s},k_{s+1}-1)
 \]
 for some $k_1,\cdots, k_{s+1}\geq 2$.
 As a result, 
 \[
 \begin{split}
 &\;\;\;\;m(\eta(M_s))       \\
 &= \sum_{k_1,\cdots,k_{s+1}\geq 2}  \left( \zeta^{\star}(k_1,\cdots,k_{s},k_{s+1}-1)- \zeta^{\star}(k_1,\cdots,k_s,k_{s+1})   \right)\\
 &=  \sum_{k_1,\cdots,k_{s+1}\geq 2}    \sum_{n_1\geq \cdots\geq n_{s}\geq n_{s+1}\geq 2} \frac{n_{s+1}-1}{n_1^{k_1}\cdots n_{s}^{k_{s}} n_{s+1}^{k_{s+1}} }  \\
 &= \sum_{n_1\geq \cdots\geq n_{s+1}\geq 2} \frac{1}{n_1(n_1-1)\cdots n_{s}(n_{s}-1) n_{s+1} } .
 \end{split}
 \]
 By Lemma \ref{leb}, we have 
 \[
      m(\eta(M_s))           \leq \frac{1}{2^{s+1}} \prod_{l\geq 3}\frac{1}{1-\frac{2}{l(l-1)}}+\zeta^\star(\{2\}^{s},1)-\zeta^\star(\{2\}^{s}).
                 \]
                 Thus 
                 \[
                  m(\eta( \mathcal{T}_2))=\mathop{\mathrm{lim}}_{s\rightarrow +\infty}    \mu(\eta(M_s)) =0.                \]

            $(ii)$ The statement $\left(\eta( \mathcal{E})\right)^\prime=[1,+\infty]$ follows immediately  from the following simple observation 
            \[
            \zeta^\star(k_1,\cdots,k_r)<\zeta^{\star}(k_1,\cdots,k_r,2,\cdots,2,\cdots)<\zeta^\star(k_1,\cdots,k_r,1 )            \]
            and Theorem \ref{rep}, Lemma \ref{plus}.
             For $r>1$, $2\leq i_1<i_2<\cdots<i_r$, define
            \[
            E_{i_1,\cdots,i_r}=\Big{\{}(k_1,\cdots,k_r,\cdots)\in\mathcal{T}\,\Big{|}\, k_{i_1}=\cdots =k_{i_r}=1, k_{i}\geq 2, i\notin \{i_1,\cdots, i_r\}\Big{\}}.
            \]
            Then 
            \[
            E_{i_1,\cdots,i_r}\cap E_{j_1,\cdots,j_s}=\emptyset, (i_1,\cdots, i_r)\neq (j_1,\cdots,j_s),
            \]
            \[
            \mathcal{E}=\bigcup_{r>1} \bigcup_{2\leq  i_1< i_2< \cdots< i_r} E_{i_1,\cdots,i_r}.
            \]
            By Theorem \ref{rep}, $x\in \eta(E_{i_1,\cdots,i_r})$ if and only if $\forall\, s>i_r$,
            \[
            x\in Z_{k_1, \cdots,k_{s-1}, k_s}
            \]
            for some  $(k_1,\cdots, k_{s-1},k_s)$ which satisfies  $$k_{i_1}=\cdots=k_{i_r}=1, k_i\geq 2, i\notin \{i_1,\cdots,i_r\}.$$
            For $s>i_r$, one  has 
            \[
            \begin{split}
            &\;\;\;\;m (\eta(E_{i_1,\cdots,i_r}))\\
            &\leq   \sum_{\substack{k_{i_1}=\cdots=k_{i_r}=1\\ k_i\geq 2, i\notin \{i_1,\cdots,i_r\}, i\leq s }}  m\left(Z_{k_1, \cdots,k_{s-1}, k_s}    \right)\\
            &\leq  \sum_{\substack{k_{i_1}=\cdots=k_{i_r}=1\\ k_i\geq 2, i\notin \{i_1,\cdots,i_r\}, i\leq s }}   \sum_{n_1\geq \cdots\geq n_s\geq 2}   \frac{n_s-1}{n_1^{k_1}\cdots n_{s-1}^{k_{s-1}}n_s^{k_s}}      \\
            &\leq \sum_{n_1\geq \cdots\geq n_s\geq 2}  \prod_{j=1}^r   \frac{1}{n_{i_j}}\cdot \prod_{\substack{ i\notin\{i_1,\cdots,i_r\}\\1\leq i<s}}\frac{1}{n_i(n_i-1)} \cdot \frac{1}{n_s}    \\
            &\leq \sum_{n_1\geq \cdots\geq n_{i_r}\geq 2}  \prod_{j=1}^r   \frac{1}{n_{i_j}}\cdot \prod_{\substack{ i\notin\{i_1,\cdots,i_r\}\\1\leq i\leq i_r}}\frac{1}{n_i(n_i-1)} \cdot  \sum_{n_{i_r}\geq n_{i_r+1}\geq \cdots \geq n_s\geq 2} \prod_{i_r<i<s}\frac{1}{n_i(n_i-1)}\cdot \frac{1}{n_s}\\
            &\leq     \sum_{n_1\geq \cdots\geq n_{i_r}\geq 2}  \prod_{j=1}^r   \frac{1}{n_{i_j}}\cdot \prod_{\substack{ i\notin\{i_1,\cdots,i_r\}\\1\leq i\leq i_r}}\frac{1}{n_i(n_i-1)} \cdot  \sum_{+\infty> n_{i_r+1}\geq \cdots \geq n_s\geq 2} \prod_{i_r<i<s}\frac{1}{n_i(n_i-1)}\cdot \frac{1}{n_s}.
   \\            \end{split}
            \]
            By Lemma \ref{leb}, it follows that
              \[
            \begin{split}
            &\;\;\;\;m(\eta(E_{i_1,\cdots,i_r}))\\
             &\leq     \sum_{n_1\geq \cdots\geq n_{i_r}\geq 2}  \prod_{j=1}^r   \frac{1}{n_{i_j}}\cdot \prod_{\substack{ i\notin\{i_1,\cdots,i_r\}\\1\leq i\leq i_r}}\frac{1}{n_i(n_i-1)} \cdot  \mathop{\mathrm{lim}}_{s\rightarrow +\infty}  \sum_{+\infty> n_{i_r+1}\geq \cdots \geq n_s\geq 2} \prod_{i_r<i<s}\frac{1}{n_i(n_i-1)}\cdot \frac{1}{n_s}.\\
             &=0.
            \end{split}
            \]
            Thus we have 
            \[
            m(\eta(E_{i_1,\cdots,i_r}))=0 .        \]
            In a word, 
            \[
             m\left(\eta(\mathcal{E})\right)=\sum_{r>1} \sum_{2\leq i_1\leq i_2\leq \cdots\leq i_r} m\left(\eta(E_{i_1,\cdots,i_r})\right) =0.           \]
             
             $(iii)$ If $x\notin \eta(\mathcal{D}_p)$, then 
 \[
 \eta^{-1}(x)=(k_1,\cdots,k_r,\cdots), k_1,\cdots, k_{i-1}\leq p, k_i\geq p+1,\;for \;some\; i\geq 1.
 \]
 By the definition of $\eta$ , we have 
 \[
 x\in{ Z_{k_1,\cdots,k_i}} =(\zeta^\star(k_1,\cdots,k_i),\zeta^\star(k_1,\cdots,k_{i-1},k_i-1)],     \]
$$Z_{k_1,\cdots,k_i}\cap \eta(\mathcal{D}_p)=\emptyset.$$   There are the following three cases:\\
$(A)$ If $x\in \left(Z_{k_1,\cdots,k_i}\right)^o$, we have $\left(Z_{k_1,\cdots,k_i}\right)^o \cap \eta(\mathcal{D}_p)=\emptyset.$\\
$(B)$ If  $x=  \zeta^\star(k_1,\cdots,k_{i-1},k_i-1)$ and $k_i\geq p+2$, then
\[x\in \left(  \zeta^\star(k_1,\cdots, k_{i-1},k_i),\zeta^\star(k_1,\cdots,k_{i-1},k_i-2)   \right)\]
and \[
\left(  \zeta^\star(k_1,\cdots, k_{i-1},k_i),\zeta^\star(k_1,\cdots,k_{i-1},k_i-2)   \right)\cap \eta(\mathcal{D}_p)=\emptyset.\]
$(C)$ If  $x=  \zeta^\star(k_1,\cdots,k_{i-1},k_i-1)$ and $k_i= p+1$, then
\[
x\in \left( \zeta^\star(k_1, \cdots,k_{i-1},p+1),\zeta^\star(k_1,\cdots, k_{i-1}, p,p)    \right)
\]
and \[
 \left( \zeta^\star(k_1, \cdots,k_{i-1},p+1),\zeta^\star(k_1,\cdots, k_{i-1}, p,p)    \right)\cap \eta(\mathcal{D}_p)=\emptyset.\]
 In a word, $\eta(\mathcal{D}_p)$ is closed.

             By the order structure of multiple zeta-star values and   \[
 \mathop{\mathrm{lim}}_{n\rightarrow +\infty} \zeta^{\star}(2, \{1\}^n)=+\infty,
 \]  
 one has 
 \[
   (1,+\infty)= \bigcup_{\substack{  k_1\geq 2,k_2,\cdots,k_r\geq 1            \\                  }              } Z_{k_1,\cdots,k_r},  \]
   \[
   (1,\zeta^\star(2,\{1\}^{t-1}]= \bigcup_{\substack{  k_1\geq 2,k_2,\cdots,k_r\geq 1            \\     (k_1,\cdots, k_t)\neq (2,\{1\}^{t-1}    )         }              } Z_{k_1,\cdots,k_r}   \]
   for $r\geq t$.
 By definition,  for $x\in \eta(\mathcal{D}_p)$, then
 \[
x\in \bigcup_{\substack{1\leq k_1,\cdots, k_r\leq p\\ k_1\geq 2          }}Z_{k_1,\cdots,k_r}
 \] for each  $r\geq 1$.
  Beware that $\eta(\mathcal{D}_p)$ is an unbounded subset of $(1,+\infty)$. For $r\geq t+2$, one has
   \[
 \begin{split}
 &\;\;\;\;m\left(\eta(\mathcal{D}_p)\cap  (1,\zeta^\star(2,\{1\}^{t-1}]   \right)\\
 &=\sum_{\substack{ 1\leq k_1,\cdots,k_r\leq p\\ k_1\geq 2, (k_1,\cdots,k_t)\neq (2,\{1\}^{t-1}      }} m\left( \eta(\mathcal{D}_p) \cap Z_{k_1,k_2,\cdots,k_r}   \right)\\
 &\leq\sum_{\substack{ 1\leq k_1,\cdots,k_r\leq p\\ k_1\geq 2, (k_1,\cdots,k_t)\neq (2,\{1\}^{t-1}      }} m\left(  Z_{k_1,k_2,\cdots,k_r}   \right)      \\
 &\leq\sum_{\substack{ 1\leq k_1,\cdots,k_r\leq p\\ k_1\geq 2, (k_1,\cdots,k_t)\neq (2,\{1\}^{t-1}      }} \sum_{n_1\geq \cdots \geq n_r\geq 2} \frac{n_r-1}{n_1^{k_1}\cdots n_r^{k_r}}   \\
  &\leq \sum_{\substack{ 1\leq k_1,\cdots,k_t\leq p\\ k_1\geq 2, (k_1,\cdots,k_t)\neq (2,\{1\}^{t-1}      }} \sum_{n_1\geq \cdots \geq n_r\geq 2}  \sum_{1\leq k_{t+1} ,\cdots, k_r\leq p}  \frac{n_r-1}{n_1^{k_1}\cdots n_r^{k_r}}   \\ 
  &\leq  \sum_{\substack{ 1\leq k_1,\cdots,k_t\leq p\\ k_1\geq 2, (k_1,\cdots,k_t)\neq (2,\{1\}^{t-1}      }} \sum_{n_1\geq \cdots \geq n_r\geq 2} \frac{1}{n_1^{k_1}\cdots n_t^{k_t}} \cdot \frac{n_r-1}{n_{t+1}\cdots n_r} \cdot \frac{(1-\frac{1}{n_{t+1}^p})\cdots   (1-\frac{1}{n_{r}^p})             }{  (1-\frac{1}{n_{t+1}})\cdots   (1-\frac{1}{n_{r}})              }\\
   &\leq  \sum_{\substack{ 1\leq k_1,\cdots,k_t\leq p\\ k_1\geq 2, (k_1,\cdots,k_t)\neq (2,\{1\}^{t-1}      }} \sum_{n_1\geq \cdots \geq n_r\geq 2} \frac{1}{n_1^{k_1}\cdots n_t^{k_t}} \cdot \frac{1}{n_{t+1}\cdots n_{r-1}} \cdot \frac{(1-\frac{1}{n_{t+1}^p})\cdots   (1-\frac{1}{n_{r}^p})             }{  (1-\frac{1}{n_{t+1}})\cdots   (1-\frac{1}{n_{r-1}})              }\\
   &\leq  \sum_{\substack{ 1\leq k_1,\cdots,k_t\leq p\\ k_1\geq 2, (k_1,\cdots,k_t)\neq (2,\{1\}^{t-1}      }} \sum_{n_1\geq \cdots \geq n_r\geq 2} \frac{1}{n_1^{k_1}\cdots n_t^{k_t}} \cdot \frac{1}{(n_{t+1}-1)\cdots (n_{r-1}-1)} \cdot \left(1-\frac{1}{n_{t+1}^p}\right)\cdots   \left(1-\frac{1}{n_{r}^p}\right)             \\     \end{split}
 \]
 By Lemma \ref{dbou}, one has 
 \[
m\left(\eta(\mathcal{D}_p)\cap  (1,\zeta^\star(2,\{1\}^{t-1}]   \right) =0\]
for all $t\geq 1$.
 In a word, for $p\geq 2$, \[
m\left(\eta(\mathcal{D}_p)    \right)=0,\]  
 \[
 m(\eta(\bigcup_{p\geq 2} \mathcal{D}_p))=0. \] 
   $\hfill\Box$\\                       
   
    \begin{lem}\label{exp}
 For $r\geq 1$, $\forall\, k_1,\cdots,k_r\geq 2$, we have 
 \[
\frac{1}{2^{k_1+k_2+\cdots+k_r}} <\sum_{n_1\geq \cdots \geq n_r\geq 2}\frac{1}{n_1^{k_1}\cdots n_{r-1}^{k_{r-1}} n_r^{k_r}}\leq \frac{C_1}{2^{k_1+k_2+\cdots+k_r}},
 \]
 \[
 \frac{1}{2^{k_1+k_2+\cdots+k_r+1}}<\sum_{n_1\geq \cdots \geq n_{r+1}\geq 2}\frac{1}{n_1^{k_1}\cdots n_{r-1}^{k_{r-1}} n_r^{k_r}n_{r+1}}\leq \frac{C_1}{2^{k_1+k_2+\cdots+k_r}}. \]
 Here the constant $C_1$ is independent of $r,k_1,\cdots,k_r$.
\end{lem}
  \noindent{\bf Proof:}  
 For $k_1,\cdots,k_{r}\geq 2$, one has 
 \[
 \begin{split}
 &\;\;\;\;  \sum_{n_1\geq \cdots \geq n_r\geq 2}\frac{1}{n_1^{k_1}\cdots n_{r-1}^{k_{r-1}} n_r^{k_r}}\\
 &=  \frac{1}{2^{k_1+\cdots+k_r}} \left(1+ \sum_{s=1}^r \sum_{n_1\geq \cdots \geq n_s\geq 3}    \left( \frac{2}{n_1}\right)^{k_1}  \cdots    \left( \frac{2}{n_s}\right)^{k_s}    \right)    \\
 &\leq  \frac{1}{2^{k_1+\cdots+k_r}} \left(1+ \sum_{s=1}^r \sum_{n_1\geq \cdots \geq n_s\geq 3}    \left( \frac{2}{n_1}\right)^{2}  \cdots    \left( \frac{2}{n_s}\right)^{2}    \right)  \\
 &\leq  \frac{1}{2^{k_1+\cdots+k_r}} \prod_{  l\geq 3}\left(1+\frac{4}{l^2}+\cdots+\left(  \frac{4}{l^2} \right)^i +\cdots           \right)\\
 &\leq  \frac{1}{2^{k_1+\cdots+k_r}} \prod_{l\geq3}\frac{1}{1-\frac{4}{l^2}} . \end{split}
  \] 
  
  By the following simple observation
  \[
  \frac{1}{l}<\frac{2}{\sqrt{l}+\sqrt{l-1}}=2\left(\sqrt{l}-\sqrt{l-1}  \right),\;\forall\, l\geq 2,
  \]
  one has 
  \[
  \sum_{l=2}^n\frac{1}{l}<2\sqrt{n}, \forall\, n\geq 2. \]
  By using the above formula, we have
   \[
   \begin{split}
&\;\;\;\; \sum_{n_1\geq \cdots \geq n_{r+1}\geq 2}\frac{1}{n_1^{k_1}\cdots n_{r-1}^{k_{r-1}} n_r^{k_r}n_{r+1}}\\
&<  \sum_{n_1\geq \cdots \geq n_{r}\geq 2}\frac{2}{n_1^{k_1}\cdots n_{r-1}^{k_{r-1}} n_r^{k_r-\frac{1}{2}}}  \\
&\leq   \frac{2}{2^{k_1+\cdots+k_r-\frac{1}{2}}} \Bigg{[}1+ \sum_{1\leq s<r} \sum_{n_1\geq \cdots \geq n_s\geq 3}    \left( \frac{2}{n_1}\right)^{k_1}  \cdots    \left( \frac{2}{n_s}\right)^{k_s}  +\\
&\;\;\;\;+ \sum_{n_1\geq \cdots \geq n_r\geq 3}    \left( \frac{2}{n_1}\right)^{k_1}  \cdots    \left( \frac{2}{n_{r-1}}\right)^{k_{r-1}}  \left( \frac{2}{n_{r}}\right)^{k_{r}-\frac{1}{2}} \Bigg{]}   \\
 &\leq  \frac{2}{2^{k_1+\cdots+k_r-\frac{1}{2}}} \left[1+ \sum_{s=1}^r \sum_{n_1\geq \cdots \geq n_s\geq 3}    \left( \frac{2}{n_1}\right)^{\frac{3}{2}}  \cdots    \left( \frac{2}{n_s}\right)^{\frac{3}{2}}    \right]  \\
 &\leq  \frac{2}{2^{k_1+\cdots+k_r-\frac{1}{2}}} \prod_{l\geq 3} \left[  1+  \left( \frac{2}{l}\right)^{\frac{3}{2}} +\cdots+  \left( \frac{2}{l}\right)^{\frac{3}{2}i}+\cdots \right] \\
 &\leq \frac{2^{\frac{3}{2}}}{2^{k_1+\cdots+k_r}} \prod _{l\geq 3} \frac{1}{1-\left( \frac{2}{l}\right)^{\frac{3}{2}}}.\\
  \end{split}
\]
Define $C_1$ as 
\[
C_1=max\Big{\{} \prod_{l\geq3}\frac{1}{1-\frac{4}{l^2}},\;\; 2^{\frac{3}{2}}    \prod _{l\geq 3} \frac{1}{1-\left( \frac{2}{l}\right)^{\frac{3}{2}}}   \Big{\}}.
\]
The lemma is proved.  $\hfill\Box$\\   

 \noindent{\bf Proof of Theorem \ref{hdim}:} $(i)$
 By definition, it is clear that
 \[
 \eta(\mathcal{T}_p)= \bigcap_{r\geq 2}\bigcup_{k_1,\cdots, k_r\geq p}Z_{k_1,\cdots,k_r}.
 \]
 For the upper bound of the Hausdorff dimension, by Lemma \ref{exp},  it suffices to show that for  $r\geq 2$  and $t=\frac{\mathrm{log}\;\alpha_p}{\mathrm{log}\;2}$,
 \[
 \sum_{k_1,\cdots, k_r\geq p}\left[m\left( Z_{k_1,\cdots,k_r}\right)\right]^t<+\infty.
  \]
  
 By Lemma \ref{exp},  for $p\geq 2$, one can check that 
  \[
  \begin{split}
  &\;\;\;\; \sum_{k_1,\cdots, k_r\geq p}\left[m\left( Z_{k_1,\cdots,k_r}\right)\right]^t\\
  &= \sum_{k_1,\cdots, k_r\geq p} \left(  \sum_{n_1\geq \cdots\geq n_r\geq 2}   \frac{n_r-1}{n_1^{k_1}\cdots n_{r-1}^{k_{r-1}}n_r^{k_r}}   \right)^t \\  \end{split}
 \]
 \[
 \begin{split}
  &<\sum_{k_1,\cdots, k_r\geq p} \left(  \sum_{n_1\geq \cdots\geq n_r\geq 2}   \frac{1}{n_1^{k_1}\cdots n_{r-1}^{k_{r-1}}n_r^{k_r-1}}   \right)^t \\
   &\leq \sum_{k_1,\cdots, k_r\geq p} \left( \frac{2C_1}{2^{k_1+\cdots+k_r}} \right)^t \\
  &\leq\frac{(2C_1)^t}{ 2^{rpt}}  \left( \sum_{k_1\geq 0}\frac{1}{2^{k_1t}}    \right) \cdots  \left( \sum_{k_r\geq 0}\frac{1}{2^{k_r t}}    \right)\\
  & \leq \frac{(2C_1)^t}{ \left[ 2^{(p-1)t}(2^t-1)\right]^r  } =(2C_1)^t. \end{split}
  \]
  Here the last equality follows from the fact that $2^t$ is the root of the equation 
  \[
  x^{p-1}(x-1)=1.
  \]
As a reuslt,
\[
 \mathrm{dim}_H\, \eta\left(\mathcal{T}_p\right)\leq  \frac{\mathrm{log}\;\alpha_p}{\mathrm{log}\;2}, \forall\,p\geq 2.
 \]
 
 For the lower bound, we construct the following map
 \[
 \beta: \mathcal{T} \rightarrow (0, \frac{1}{2})
 \]
 by 
 \[
 \beta\left((k_1,\cdots, k_r,\cdots)     \right)=\frac{1}{2^{k_1}}+\frac{1}{2^{k_1+k_2}}+\cdots+\frac{1}{2^{k_1+\cdots+k_r}}+\cdots.
 \]
 By the binary expansion of the real number, the map $\beta$ is bijective.
 What is more, for any ${\bf k},{\bf m}\in \mathcal{T}$, 
 \[
{\bf k} \succ {\bf m}
 \]
 if and only if 
 \[
 \beta({ \bf k}) >\beta ( {\bf m}). 
 \]
 Since the maps $\eta$ and $\beta$ are both bijective, one can define the map $\tau=\beta\circ \eta^{-1}$ as 
 \[
 \tau: (1,+\infty)\rightarrow (0,\frac{1}{2}),
 \]
 \[
\mathop{\mathrm{lim}}_{r\rightarrow +\infty}\zeta^{\star}(k_1,\cdots,k_r) \mapsto \frac{1}{2^{k_1}}+\cdots+\frac{1}{2^{k_1+\cdots+k_r}}+\cdots.
  \] 
  As the maps $\eta$ and $\beta$ both preserve the order structure, it follows that the map $\tau$ is a homeomorphism  of topology spaces. 
    For $p\geq 2$ and $x,y\in \eta(\mathcal{T}_p)$, by Lemma \ref{exp}, one has 
    \[
    c_2 |x-y|\leq |\tau(x)-\tau(y)|\leq c_3 |x-y|.
    \]
  Here the constants $c_2, c_3$ are independent of $p,x,y$. By Proposition {3.3} in \cite{fal},
  \[
  \mathrm{dim}_H \eta(\mathcal{T}_p)= \mathrm{dim}_H \tau\left(\eta(\mathcal{T}_p)\right) = \mathrm{dim}_H \beta(\mathcal{T}_p).\]
  
  By the Chapter $3.4$ in \cite{fal},   a  binary interval is of the form 
  \[
  \left[\frac{i}{ 2^{k}}, \frac{i+1}{ 2^{k}}\right]
  \]
  where $r=0,1,2,\cdots , 2^{k-1}, k\geq 1$. From the discussion of Chapter $3.4$ in \cite{fal}, one can give both the upper and lower bound of   the  $t$-dimensional Hausdorff measure  $ \mathcal{H}^t\left( \beta(\mathcal{T}_p) \right)$ of $\beta(\mathcal{T}_p)$ by using the $\delta$-covers of $\beta(\mathcal{T}_p)$ by binary intervals. 
  
  By definition, $$\beta(\mathcal{T}_p)=\bigcap_{r\geq 2}\bigcup_{k_1,\cdots,k_r\geq p} \mathcal{U}_{k_1,\cdots, k_r},$$
  where \[
  \mathcal{U}_{k_1,\cdots,k_r}=\left(\frac{1}{2^{k_1}}+\cdots+\frac{1}{2^{k_1+\cdots +k_r}},   \frac{1}{2^{k_1}}+\cdots+\frac{1}{2^{k_1+\cdots +k_{r-1}}}+ \frac{1}{2^{k_1+\cdots +k_{r}-1}}       \right).
  \]
  By the theory of binary expansion of real numbers and Chapter $3.4$, \cite{fal}, we have 
   \[
  \mathcal{H}^t\left( \beta(\mathcal{T}_p) \right)\geq \frac{1}{2^{t+1}} \mathop{\mathrm{lim}}_{r\rightarrow +\infty} \sum_{k_1,\cdots,k_r\geq p}\left( \mu(\mathcal{U}_{k_1,\cdots,k_r})   \right)^t=\frac{1}{2^{t+1}}.
  \]
  Thus \[
 \mathrm{dim}_H\, \beta\left(\mathcal{T}_p\right)\geq  \frac{\mathrm{log}\;\alpha_p}{\mathrm{log}\;2}, \forall\,p\geq 2.
 \]
 In a word, we have 
 \[
 \mathrm{dim}_H\, \eta\left(\mathcal{T}_p\right)=  \frac{\mathrm{log}\;\alpha_p}{\mathrm{log}\;2}, \forall\,p\geq 2.
 \]
 
 $(ii)$ Denote by $u= \frac{\mathrm{log}\;\frac{1}{\gamma_{p,q}}}{\mathrm{log}\;2}$.  For $p<q$, as the map $\tau$ is bi-Lipschitz on the set $\eta(\mathcal{T}_p)$, one has 
 \[
  \mathrm{dim}_H\, \eta\left(\mathcal{T}_p \cap \mathcal{D}_q\right)=  \mathrm{dim}_H\,\tau\left( \eta\left(\mathcal{T}_p \cap \mathcal{D}_q\right) \right)= \mathrm{dim}_H\, \beta\left(\mathcal{T}_p \cap \mathcal{D}_q\right).  \]
    It suffices to show that 
   \[
  \mathrm{dim}_H\beta\left(  \mathcal{T}_p \cap \mathcal{D}_q\right)=u.
    \] 
    It is clear that
    $$\beta(   \mathcal{T}_p \cap \mathcal{D}_q)=\bigcap_{r\geq 2}\bigcup_{\substack{k_1\geq 2 \\p\leq k_1,\cdots,k_r\leq q }} \mathcal{U}_{k_1,\cdots, k_r}.$$
    Here
    \[
    \begin{split}
  &\;\;\;\mathcal{U}_{k_1,\cdots,k_r}\\
  &=\begin{cases}
  \left(\frac{1}{2^{k_1}}+\cdots+\frac{1}{2^{k_1+\cdots +k_r}},   \frac{1}{2^{k_1}}+\cdots+\frac{1}{2^{k_1+\cdots +k_{r-1}}}+ \frac{1}{2^{k_1+\cdots +k_{r}-1}}       \right], & k_r\geq 2;\\
  \left(\frac{1}{2^{k_1}}+\cdots+\frac{1}{2^{k_1+\cdots +k_r}},   \frac{1}{2^{k_1}}+\cdots+\frac{1}{2^{k_1+\cdots +k_{i-1}}}+ \frac{1}{2^{k_1+\cdots +k_{i-1}+k_i-1}}       \right],&k_i\geq 2, k_{i+1}=\cdots=k_r=1.\\
  \end{cases}.\\
  \end{split}
  \]
   For any $r\geq 2$, one has 
   \[
   \begin{split}
   &\;\;\;\;\sum_{p\leq k_1,\cdots,k_r\leq q }\left[\mu( \mathcal{U}_{k_1,\cdots, k_r})\right]^u \\
           &=   \sum_{p\leq k_1,\cdots,k_r\leq q }\frac{1}{2^{(k_1+\cdots+k_r)u}}\\
         &=\left( \sum_{p\leq  k_1 \leq q}\frac{1}{2^{ k_1 u}}\right)^r  \\
         &=1.
    \end{split}
      \]
   Here the last equality follows from the fact that $\frac{1}{2^{u}}$ is the root of the equation 
   \[
   x^p+x^{p+1}+\cdots+x^q=1.
   \]
   Since $\mu( \mathcal{U}_{k_1,\cdots, k_r})=\frac{1}{2^{k_1+\cdots+k_r}}\rightarrow 0$ as $r\rightarrow +\infty$, \[
  \mathrm{dim}_H\beta\left(  \mathcal{T}_p \cap \mathcal{D}_q\right)\leq u.
    \] 
     By the theory of binary expansion of real numbers and Chapter $3.4$, \cite{fal}, we have 
   \[
  \mathcal{H}^u\left( \beta(\mathcal{T}_p\cap \mathcal{D}_q) \right)\geq \frac{1}{2^{u+1}} \mathop{\mathrm{lim}}_{r\rightarrow +\infty} \sum_{q\geq k_1,\cdots,k_r\geq p}\left( \mu(\mathcal{U}_{k_1,\cdots,k_r})   \right)^u=\frac{1}{2^{u+1}}.
  \]
  So 
  \[
    \mathrm{dim}_H\beta\left(  \mathcal{T}_p \cap \mathcal{D}_q\right)\geq u.  \]
    As a result, 
    \[
      \mathrm{dim}_H\eta\left(  \mathcal{T}_p \cap \mathcal{D}_q\right)= u .   \]
      
      $(iii)$ The statement $(iii)$ follows immediately from $(ii)$.
   
      $\hfill\Box$\\       
  
  \begin{rem}
  For $p\geq 2$, $\eta(\mathcal{T}_p)$ is a bounded closed set of measure zero, while $\eta(\mathcal{D}_p)$ is an unbounded set of measure zero. Since 
  \[
  \eta(\mathcal{T}_{p+1})\cap \eta(\mathcal{D}_p)=\emptyset, \forall\,p\geq 2,    \]
  for $p\geq 2$, what is the structure of 
  \[
   \eta(\mathcal{T}_{p+1})+ \eta(\mathcal{D}_p) ? \]
   Similarly, do the sets
   \[
    \eta(\mathcal{T}_{p})+ \eta(\mathcal{T}_q),   \eta(\mathcal{D}_{p})+ \eta(\mathcal{D}_q),  \eta(\mathcal{T}_{p})+ \eta(\mathcal{D}_q) \]
   have interior points?
    \end{rem}       
    
    \begin{rem}
    By using the binary intervals, one can compute
    \[
    \mathrm{dim}_H \beta (\mathcal{D}_2).
    \]
    Since the map $\tau$ is not a bi-Lipschitz transformation on the set $\eta(\mathcal{D}_2)$, we don't  knowthe  Hausdorff dimension of $\eta(\mathcal{D}_2)$. One can prove that 
    \[
     \mathrm{dim}_H \eta (\mathcal{D}_2)\geq   \mathrm{dim}_H \beta (\mathcal{D}_2) .  
       \]
      It seems that    $   \mathrm{dim}_H \eta (\mathcal{D}_2)>  \mathrm{dim}_H \beta (\mathcal{D}_2) $.    \end{rem}               
    
    \begin{rem}
    For $p=2$, then $\alpha_p=\frac{\sqrt{5}+1}{2}$
    and 
   \[
 \mathrm{dim}_H\, \eta\left(\mathcal{T}_2\right)= \frac{\mathrm{log}\; \frac{\sqrt{5}+1}{2}}{\mathrm{log}\;2}. \]
 This is equal to the Hausdorff dimension of the following  Cantor set
 \[
 C_{[\frac{1}{4},\frac{2}{4}]}=\Big{\{} x\in [0,1]\,\Big{|}\,x=\frac{i_1}{4}+\frac{i_2}{4}+\cdots+\frac{i_r}{4^r}+\cdots , i_j \in \{0,2,3\}    \Big{ \}}.
 \]
    \end{rem}    
  
\section{Another approach to the multiple integral representations}\label{sie}
 
 In this section, we will give another approach to the multiple integrals in Theorem \ref{last}.
 
 From Theorem \ref{last}, it follows that
 \[
 \mathop{\int}_{[0,1]^{2k}}\frac{dx_1dx_2\cdots dx_{2k}}{1-x_1+\cdots +(-1)^ix_1\cdots x_i+\cdots+x_1\cdots x_{2k}}=\zeta^{\star}(\underbrace{2,\cdots,2}_{k})
, \tag{7}\]
\[
\mathop{\int}_{[0,1]^{2k+1}}\frac{dx_1dx_2\cdots dx_{2k+1}}{1-x_1+\cdots +(-1)^ix_1\cdots x_i+\cdots-x_1\cdots x_{2k+1}}=\zeta^{\star}(\underbrace{2,\cdots,2}_{k},1 )
.\tag{8}\]
By the results of \cite{oz}, we have
\[
\zeta^{\star}(\underbrace{2,\cdots,2}_{k})=2(1-2^{1-2k})\zeta(2k),\; \zeta^{\star}(\underbrace{2,\cdots,2}_{k},1)=2\zeta(2k+1). \tag{9}
\]
Thus the multiple integrals $(7)$ and $(8)$ are in fact Riemann zeta values.

 Now we give an alternative way to calculate the multiple integrals $(5)$ and $(6)$.
\begin{prop}\label{zeta}
For $k\geq 1$, \\
$(i)$ \[
\zeta(2k)=\frac{1}{2(1-2^{1-2k})}\mathop{\int}_{\Delta_{2k}}\frac{dt_1}{1-t_1}\frac{dt_2}{t_2}\cdots \frac{dt_{2k-1}}{1-t_{2k-1}}\frac{dt_{2k}}{t_{2k}}
,\]
where 
\[
\Delta_{2k}=\{(t_1,\cdots,t_{2k})\in (0,1)^{2k} \;\big{|} \;t_1<t_2,t_{2i-1}<t_{2i-2},t_{2i-1}<t_{2i}, i=2,\cdots,k\};
\]
$(ii)$ 
\[
\zeta(2k+1)=\frac{1}{2}\mathop{\int}_{\Delta_{2k+1}}\frac{dt_1}{1-t_1}\frac{dt_2}{t_2}\cdots \frac{dt_{2k-1}}{1-t_{2k-1}}\frac{dt_{2k}}{t_{2k}}\frac{dt_{2k+1}}{1-t_{2k+1}},\]
where $$\Delta_{2k+1}=\{(t_1,\cdots,t_{2k+1})\in (0,1)^{2k+1}\; \big{|}\; t_{2i}>t_{2i-1}, t_{2i}>t_{2i+1}, i=1,\cdots,k\}.$$
\end{prop}
   \noindent{\bf Proof:} 
   $(i)$ For $(x_1,\cdots,x_{2k})\in (0,1)^{2k}$, by changing of variables:
  \[t_1=x_1-x_1x_2+\cdots +x_1x_2\cdots x_{2k-1}-x_1x_2\cdots x_{2k},          \]
  \[
  t_2=1-x_2+\cdots +x_2\cdots x_{2k-1}-x_2\cdots x_{2k},
  \]
  \[
  \cdots\;\;\;\;\;\;\cdots\;\;\;\;\;\;\cdots
  \]
  \[
  t_{2k-1}=x_{2k-1}-x_{2k-1}x_{2k},
  \]
  \[
  t_{2k}=1-x_{2k},
  \]
  it is easy to check that:
  \[
  t_1<t_2>t_3<t_4>\cdots>t_{2k-1}<t_{2k},\; (t_1,\cdots, t_{2k})\in (0,1)^{2k}
  \]
  and 
  \[
  \frac{dx_1dx_2\cdots dx_{2k}}{1-x_1+x_1x_2+\cdots +x_1x_2\cdots x_{2k}}=\frac{dt_1}{1-t_1}\frac{dt_2}{t_2}\cdots \frac{dt_{2k-1}}{1-t_{2k-1}}\frac{dt_{2k}}{t_{2k}}.
  \]
  Thus 
  \[
   \mathop{\int}_{[0,1]^{2k}}\frac{dx_1dx_2\cdots dx_{2k}}{1-x_1+x_1x_2+\cdots +x_1x_2\cdots x_{2k}}=\mathop{\int}_{\Delta_{2k}}\frac{dt_1}{1-t_1}\frac{dt_2}{t_2}\cdots \frac{dt_{2k-1}}{1-t_{2k-1}}\frac{dt_{2k}}{t_{2k}}.
     \]
     
     As a result, the statement $(i)$ follows from $(5)$ and $(7)$.
   Similarly, by changing of variables:
     \[t_1=x_1-x_1x_2+\cdots +x_1x_2\cdots x_{2k-1}-x_1x_2\cdots x_{2k}+x_1x_2\cdots x_{2k+1},          \]
  \[
  t_2=1-x_2+\cdots +x_2\cdots x_{2k-1}-x_2\cdots x_{2k}+x_2\cdots x_{2k+1},
  \]
  \[
  \cdots\;\;\;\;\;\;\cdots\;\;\;\;\;\;\cdots
  \]
  \[
  t_{2k-1}=x_{2k-1}-x_{2k-1}x_{2k}+x_{2k-1}x_{2k}x_{2k+1},
  \]
  \[
  t_{2k}=1-x_{2k}+x_{2k}x_{2k+1},
  \]
  \[
  t_{2k+1}=x_{2k+1},
  \]
   one can prove the statement $(ii)$.
     $\hfill\Box$\\
     
     The simple observation in Proposition \ref{zeta} can be applied to more general cases.
     
     For $r=2k+1$, $1\leq i_1<\cdots<i_r$,  from Theorem \ref{last}, one has
\[
\begin{split}
&\;\;\mathop{\int}_{[0,1]^{i_r}}\frac{dx_1\cdots dx_{i_r}}{1-x_1\cdots x_{i_1}+x_1\cdots x_{i_1}x_{i_1+1}\cdots x_{i_2}+\cdots+(-1)^rx_1\cdots x_{i_r}}\\
&= \zeta^{\star}( i_1+1   , \underbrace{1,\cdots,1}_{i_2-i_1-1},    i_3-i_2+1, \underbrace{1,\cdots, 1}_{i_4-i_3-1},  \cdots,   i_{2k-1}-i_{2k-2}+1,  \underbrace{1,\cdots,1}_{i_{2k}-i_{2k-1}-1},  i_{2k+1}-i_{2k}      ).
\end{split}
\]
Let 
\[
t_1=x_1\cdots x_{i_1}-x_1\cdots x_{i_2}+\cdots +x_1\cdots x_{i_{2k+1}},
\]
\[
\cdots\;\;\;\;\;\;\;\;\;\;\;\;\;\cdots \;\;\;\;\;\;\;\;\;\;\;\;\;\cdots
\]
\[
t_{i_1}=x_{i_1}-x_{i_1}\cdots x_{i_2}+\cdots +x_{i_1}\cdots x_{i_{2k+1}},
\]
\[
t_{i_1+1}=1-x_{i_1+1}\cdots x_{i_2}+\cdots +x_{i_1+1}\cdots x_{i_{2k+1}},
\]
\[
\cdots\;\;\;\;\;\;\;\;\;\;\;\;\;\cdots \;\;\;\;\;\;\;\;\;\;\;\;\;\cdots
\]
\[
t_{i_2}=1-x_{i_2}+\cdots +x_{i_2}\cdots x_{i_{2k+1}},
\]
\[
t_{i_2+1}=x_{i_2+1}\cdots x_{i_3}-x_{i_2+1}\cdots x_{i_4}+\cdots +x_{i_2+1}\cdots x_{i_{2k+1}},
\]
\[
\cdots\;\;\;\;\;\;\;\;\;\;\;\;\;\cdots \;\;\;\;\;\;\;\;\;\;\;\;\;\cdots
\]
\[
t_{i_3}=x_{i_3}-x_{i_3}\cdots x_{i_4}+\cdots +x_{i_3}\cdots x_{i_{2k+1}},
\]
\[
\cdots\;\;\;\;\;\;\;\;\;\;\;\;\;\cdots \;\;\;\;\;\;\;\;\;\;\;\;\;\cdots
\]
\[
t_{i_{2k}+1}=x_{i_{2k}+1}\cdots x_{i_{2k+1}},
\]
\[
\;\;\;\;\;\;\;\;\;\;\;\;\;\cdots \;\;\;\;\;\;\;\;\;\;\;\;\;
\]
\[
t_{i_{2k+1}}=x_{i_{2k+1}}.
\]
For $(x_1,\cdots,x_{i_{2k+1}})\in (0,1)^{i_{2k+1}}$, it is easy to check that 
\[
t_1<\cdots <t_{i_1}<t_{i_1+1}>\cdots>t_{i_2}>t_{i_2+1}<\cdots>t_{i_{2k}+1}<\cdots<t_{i_{2k+1}}
\]
and  $$(t_1,\cdots,t_{i_{2k+1}})\in (0,1)^{i_{2k+1}}$$
 Moreover,
 \[
 \begin{split}
 \frac{dt_1}{t_2}\cdots \frac{dt_{i_1}}{t_{i_1+1}}\frac{dt_{i_1+1}}{1-t_{i_1+2}}\cdots \frac{dt_{i_2}}{1-t_{i_2+1}}\cdots \frac{dt_{i_{2k}+1}}{t_{i_{2k}+2}}\cdots \frac{dt_{i_{2k+1}-1}}{t_{i_{2k+1}}}\cdot dt_{i_{2k+1}}=dx_1dx_2\cdots dx_{i_{2k+1}}.
 \end{split}
 \]
 Thus we have
 \[
\begin{split}
&\;\;\zeta^{\star}(i_1+1, \underbrace{1,\cdots,1}_{i_2-i_1-1}, i_3-i_2+1, \underbrace{1,\cdots, 1}_{i_4-i_3-1},\cdots, \underbrace{1,\cdots,1}_{i_{2k}-i_{2k-1}-1}, i_{2k+1}-i_{2k}   )              \\
&=\mathop{\int}_{[0,1]^{i_{2k+1}}}\frac{dx_1\cdots dx_{i_{2k+1}}}{1-x_1\cdots x_{i_1}+x_1\cdots x_{i_1}x_{i_1+1}\cdots x_{i_2}+\cdots-x_1\cdots x_{i_{2k+1}}}\\
&=\mathop{\int}_{\Omega_{i_1,i_2,\cdots,i_{2k+1}}}\frac{dt_1}{1-t_1}\frac{dt_2}{t_2}\cdots \frac{dt_{i_1}}{t_{i_1}}\frac{dt_{i_1+1}}{t_{i_1+1}}\frac{dt_{i_1+2}}{1-t_{i_1+2}}\cdots \frac{dt_{i_2}}{1-t_{i_2}}\cdots \frac{dt_{i_{2k}+1}}{1-t_{i_{2k}+1}}\frac{dt_{i_{2k}+2}}{t_{i_{2k}+2}}\cdots \frac{dt_{i_{2k+1}}}{t_{i_{2k+1}}},
\end{split}
\]
where 
\[
\begin{split}
&\Omega_{i_1,i_2,\cdots, i_{2k+1}}=\Big{\{}(t_1,t_2,\cdots,t_{i_{2k+1}})\in (0,1)^{i_{2k+1}}\;\\
&\;\;\;\;\;\;\;\;\;\;\;\;\big{|}\;
 t_1<\cdots <t_{i_1}<t_{i_1+1}>\cdots>t_{i_2}>t_{i_2+1}<\cdots>t_{i_{2k}+1}<\cdots<t_{i_{2k+1}} \Big{\}}.\\
\end{split}
\]
By the same way, one can deduce that
 \[
\begin{split}
&\;\;\zeta^{\star}( i_1+1   , \underbrace{1,\cdots,1}_{i_2-i_1-1},    i_3-i_2+1, \underbrace{1,\cdots, 1}_{i_4-i_3-1},  \cdots,   i_{2k-1}-i_{2k-2}+1,  \underbrace{1,\cdots,1}_{i_{2k}-i_{2k-1}-1},  i_{2k+1}-i_{2k}      ) \\
&=\mathop{\int}_{[0,1]^{i_{2k}}}\frac{dx_1\cdots dx_{i_{2k}}}{1-x_1\cdots x_{i_1}+x_1\cdots x_{i_1}x_{i_1+1}\cdots x_{i_2}+\cdots-x_1\cdots x_{i_{2k}}}\\
&=\mathop{\int}_{\Omega_{i_1,i_2,\cdots,i_{2k}}}\frac{dt_1}{1-t_1}\frac{dt_2}{t_2}\cdots \frac{dt_{i_1}}{t_{i_1}}\frac{dt_{i_1+1}}{t_{i_1+1}}\frac{dt_{i_1+2}}{1-t_{i_1+2}}\cdots \frac{dt_{i_2}}{1-t_{i_2}}\cdots \frac{dt_{i_{2k-1}+1}}{t_{i_{2k-1}+1}}\frac{dt_{i_{2k-1}+2}}{1-t_{i_{2k-1}+2}}\cdots \frac{dt_{i_{2k}}}{1-t_{i_{2k}}},
\end{split}
\]
where 
\[
\begin{split}
&\Omega_{i_1,i_2,\cdots, i_{2k}}=\Big{\{}(t_1,t_2,\cdots,t_{i_{2k}})\in (0,1)^{i_{2k}}\;\\
&\;\;\;\;\;\;\;\;\;\;\;\;\big{|}\;
 t_1<\cdots <t_{i_1}<t_{i_1+1}>\cdots>t_{i_2}>t_{i_2+1}<\cdots<t_{i_{2k-1}+1}>\cdots>t_{i_{2k}} \Big{\}}.\\
\end{split}
\]
Denote by $S_{2k}, S_{2k+1}$  the permutation group of $\{1,\cdots,2k\}$ and $\{1,\cdots, 2k+1\}$ respectively.
We have
\[
\Omega_{i_1,i_2,\cdots,i_{2k}}=\bigcup_{\sigma\in\Sigma_{2k}\subseteq S_{2k}}\{(t_1,\cdots,t_{2k})\;|\;0<t_{\sigma(1)}<t_{\sigma(2)}<\cdots<t_{\sigma(2k)}<1\} \bigcup R_{2k},
\]
\[
\Omega_{i_1,i_2,\cdots,i_{2k+1}}=\bigcup_{\sigma\in\Sigma_{2k+1}\subseteq S_{2k+1}}\{(t_1,\cdots,t_{2k+1})\;|\;0<t_{\sigma(1)}<t_{\sigma(2)}<\cdots<t_{\sigma(2k+1)}<1\} \bigcup R_{2k+1}.
\]
 Here 
 \[
 \mathrm{dim}\;R_{2k}<2k,\;\mathrm{dim}\; R_{2k+1}<2k+1.
 \]
and $\Sigma_{2k}$ and $\Sigma_{2k+1}$ are subsets of $S_{2k}$ and $S_{2k+1}$ respectively. 
Through the above analysis, we essentially give a different approach to Yamamoto's Theorem $1.2$ in \cite{yam}.
    
    \section{Calculations of some special cases}
    In this section we will calculate the images of $\eta$ and $\xi$ for some special elements. 
    
    \begin{prop}
    For $k_1,p\geq 2, k_2,\cdots,k_r\geq 1$, 
    \[
    \mathop{\mathrm{lim}}_{m\rightarrow +\infty}    \zeta^{\star}(k_1,\cdots,k_{r}, \{p\}^m)=\zeta^{\star}(k_1,\cdots,k_{r-1})+ \sum_{n_1\geq\cdots\geq n_r\geq 2}   \frac{1}{n_1^{k_1}\cdots n_r^{k_r}} \prod_{n_r\geq l\geq 2}\frac{l^p}{l^p-1} .       \]
    \end{prop}
     \noindent{\bf Proof:}  By Lemma \ref{min},
     \[       
     \begin{split}
     &\;\;\;\;  \mathop{\mathrm{lim}}_{m\rightarrow +\infty}    \zeta^{\star}(k_1,\cdots,k_{r}, \{p\}^m)             \\
     &=   \zeta^{\star}(k_1,\cdots,k_{r})+\sum_{m=1}^{+\infty}  \left( \zeta^{\star}(k_1,\cdots,k_{r}, \{p\}^m)  - \zeta^{\star}(k_1,\cdots,k_{r}, \{p\}^{m-1})   \right)               \\
     &=  \zeta^{\star}(k_1,\cdots,k_{r})+ \sum_{m=1}^{+\infty} \sum_{n_1\geq\cdots\geq n_r\geq n_{r+1}\geq \cdots\geq n_{r+m}\geq 2}\frac{1}{n_1^{k_1}\cdots n_r^{k_r}n_{r+1}^p\cdots n_{r+m}^p } \\
     &= \zeta^{\star}(k_1,\cdots,k_{r-1})+ \sum_{n_1\geq\cdots\geq n_r\geq 2}   \frac{1}{n_1^{k_1}\cdots n_r^{k_r}} \left(1+ \sum_{m=1}^{+\infty} \sum_{n_r\geq  n_{r+1}\geq \cdots\geq n_{r+m}\geq 2} \frac{1}{n_{r+1}^p\cdots n_{r+m}^p }  \right)    \\
     &=  \zeta^{\star}(k_1,\cdots,k_{r-1}) +\sum_{n_1\geq\cdots\geq n_r\geq 2}   \frac{1}{n_1^{k_1}\cdots n_r^{k_r}} \prod_{n_r\geq l\geq 2}\left( \sum_{j\geq 0}\frac{1}{l^{jp}}    \right)    \\
     &=\zeta^{\star}(k_1,\cdots,k_{r-1})+ \sum_{n_1\geq\cdots\geq n_r\geq 2}   \frac{1}{n_1^{k_1}\cdots n_r^{k_r}} \prod_{n_r\geq l\geq 2}\frac{l^p}{l^p-1} . \end{split}
     \]
      $\hfill\Box$\\     
      \begin{Cor}
      $(i)$  \[
       \mathop{\mathrm{lim}}_{m\rightarrow +\infty}    \zeta^{\star}(k_1,\cdots,k_{r}, \{2\}^m)=  \sum_{n_1\geq\cdots\geq n_r\geq 1}\frac{2}{n_1^{k_1}\cdots n_{r-1}^{k_{r-1}}n_r^{k_r-1}(n_r+1)}  ;   \]
      $(ii)$   \[\mathop{\mathrm{lim}}_{m\rightarrow +\infty}    \zeta^{\star}( \{p\}^m)=    \prod_{n\geq 2} \frac{1}{1-\frac{1}{n^p}}.\]  \end{Cor} 
      
      From \cite{OW}, for $l,m\geq 1$, one has 
      \[
      \zeta^\star(\{2,\{1\}^{m-1}\}^{l},1)=(m+1)\zeta(l(m+1)+1).
      \]
      By Lemma \ref{plus} and Theorem \ref{rep}, $(i)$, it follows that
      \[
        \mathop{\mathrm{lim}}_{l \rightarrow +\infty} \zeta^\star(\{2,\{1\}^{m-1}\}^{l})=m+1.      \]
        By Theorem \ref{rep}, $(ii)$ and Lemma \ref{sep}, $(ii)$, we have the following result:
         \begin{Thm}
 Every multiple zeta-star value is not an integer.
 \end{Thm}
 
 \begin{rem}
 For $k_r\geq r+1, \forall\, r\geq 2$, one has 
 \[
\Big{|} \zeta^\star(k_1,\cdots,k_r)-\zeta^\star(k_1,\cdots,k_{r-1})\Big{|}=\sum_{n_1\geq \cdots n_r\geq 2}\frac{1}{n_1^{k_1}\cdots n_r^{k_r}}\leq \frac{C_1}{2^{\frac{r(r+3)}{2}}}.
 \]
 Thus the limit 
 \[
 \eta((k_1,\cdots,k_r,\cdots))=\mathop{\mathrm{lim}}_{r\rightarrow +\infty} \zeta^\star(k_1,\cdots,k_r) \]
 converges  fastly. By the theory of Liouville numbers, $\beta((k_1,\cdots,k_r,\cdots))$ is transcendental.
 Is the number $ \eta((k_1,\cdots,k_r,\cdots)) $ irrational or transcendental? As each $m+1$ is the image of a bounded sequence.
 \end{rem}

 \section*{Acknowledgements}
         The author wants to thank Shengyou Wen and Yufeng Wu for helpful discussions about fractal geometry.  The author thanks Yasuo Ohno and Yuta Kadono    for the notification of Zlobin's paper.     This project is  supported   the National Natural Science Foundation of China (Grant No. 12571009).

\end{document}